\documentclass[10pt]{article}
\usepackage{amsmath,amssymb,amsthm,mathrsfs,dsfont}

\usepackage[margin=3cm]{geometry} 
\usepackage{titlesec,hyperref}

\usepackage{color}

\usepackage{fancyhdr}
\titleformat{\subsection}{\it}{\thesubsection.\enspace}{1pt}{}

\newtheorem{theo}{Theorem}[section]
\newtheorem{lemm}[theo]{Lemma}

\newtheorem{coro}[theo]{Corollary}

\newtheorem{rema}[theo]{Remark}
\numberwithin{equation}{section}

\allowdisplaybreaks 

\newcommand\lm{{\lesssim}}

\begin{document}
\title{Global Well-Posedness for the 3-D Inviscid Oldroyd-B Model with Fractional Stress Diffusion
\hspace{-4mm}
}

\author{Yi $\mbox{Lu}^1$ \footnote{Email: luyi76@mail2.sysu.edu.cn},\quad
	Chen $\mbox{Liang}^1$\footnote{E-mail: liangchen89@mail2.sysu.edu.cn  } \quad and\quad
	Zhaoyang $\mbox{Yin}^{1}$\footnote{E-mail: mcsyzy@mail.sysu.edu.cn}\\
	$^1\mbox{School}$ of Science,\\ Shenzhen Campus of Sun Yat-sen University, Shenzhen 518107, China}
\date{}
\maketitle
\hrule

\begin{abstract}
We study the 3-D incompressible inviscid Oldroyd-B model with fractional stress diffusion \((-\Delta)^\beta\tau\) and without stress damping. For $\beta
\in (\frac{1}{2},1)$ and $s\in [2\beta+\frac{1}{2},\frac{5}{2}]$, we establish local well-posedness for initial data $(u_0,\tau_0)$ in \(H^s(\mathbb R^3)\), provided that \(\|(u_0,\tau_0)\|_{H^\beta}\) is bounded and \(\|\Lambda^\beta(u_0,\tau_0)\|_{H^{s-\beta}}\) is sufficiently small. If, in addition, $\beta\in(\frac{1}{2},\frac{5}{6}]$, we give a global well-posedness result. Our results complement the work of Constantin, Wu, Zhao and Zhu
\cite{ConstantinWuZhaoZhu2021} by treating the Sobolev regularity
$s\in[2\beta+\frac12,\frac{5}{2}]$, below the threshold $s>\frac52$ considered
there. We also provide an $\ell^2$ alternative to the
high-frequency $\ell^1$ summability condition in the framework
of Wu and Zhao \cite{WuZhao2022}.

\vspace*{5pt}
\noindent {\it 2020 Mathematics Subject Classification}: 35Q31, 76A05, 74B20, 42A38.

\vspace*{5pt}
\noindent{\it Keywords}: The inviscid Oldroyd-B model; Fractional stress tensor diffusion; Global well-posedness.
\end{abstract}

\vspace*{-20pt}

	\section{Introduction}

 In this paper, we study the incompressible Oldroyd-B model of the non-Newtonian fluid in \(\mathbb{R}^+ \times \mathbb{R}^d\),
	\begin{align}\label{eq1}
		\left\{\begin{array}{l}
			\partial_tu + u\cdot\nabla u-\nu\Delta u+\nabla {\rm P} = \mu_1{\rm div}~\tau,~~~~{\rm div}~u=0,\\[1ex]
			\partial_t\tau + u\cdot\nabla\tau+a\tau+Q(\nabla u,\tau)-\eta\Delta\tau=\mu_2 D(u),\\[1ex]
			u|_{t=0}=u_0,~~\tau|_{t=0}=\tau_0, \\[1ex]
		\end{array}\right.
	\end{align}
	where \(u\) denotes the velocity, \(\tau = \tau_{i,j}\) is the non-Newtonian part of the stress tensor(\(\tau\) is a \(d \times d\) symmetric matrix here) and \(\rm P\) is a scalar pressure of fluid. \(D(u)\) is the symmetric part of \(\nabla u\), that is,
\[
D(u) = \frac{1}{2} (\nabla u + (\nabla u)^T).
\]
The bilinear term \(Q\) has the following form:
\[
Q(\tau, \nabla u) = \tau \Omega(u) - \Omega(u)\tau + b(\mathrm{D}(u)\tau + \tau \mathrm{D}(u)),
\]
where \(b\) is a parameter in \([-1, 1]\), \(\Omega(u)\) is the skew-symmetric part of \(\nabla u\), i.e.,
\[
\Omega(u) = \frac{1}{2} (\nabla u - (\nabla u)^T).
\]
The Oldroyd-B model is used to describe the motion of some viscoelastic flows. In the pioneer paper \cite{Oldroyd1950,Oldroyd1958}, Oldroyd first introduced the formulations about viscoelastic flows of Oldroyd-B type. For more detailed physical background and derivations about this model, we bring the reader's attention to the papers \cite{HuLelievre2007,Oldroyd1950,Oldroyd1958}.

The well-posedness of the system (1.1) has been widely investigated in the references. In the papers \cite{GuillopeSaut1990,GUILLOPE1990849}, local solutions with large initial data were obtained by Guillopé and Saut. Moreover, it was shown there that these solutions are global when the coupling and the initial data are small enough. In the corotational case, i.e. \(b = 0\), Lions and Masmoudi established the existence of global weak solution in \cite{LionsMasmoudi2000}. However, the case \(b \neq 0\) is still an open problem by now. In the framework of the critical Besov spaces, Chemin and Masmoudi \cite{CheminMasmoudi2001} first obtained the local solutions when \(\nu > 0\), \(\mu_1 > 0\), \(\mu_2 > 0\) and \(\eta = 0\) and showed these solutions are global provided that the initial data are small enough and \(\mu_1 \mu_2 \leq c a \nu\), which implies that the coupling effect between the two equations is less important than the effect of viscosity. Later, these results were improved by Zi, Fang and Zhang \cite{ZiFangZhang2014} in the critical \(L^p\) framework, in particular, the restriction that the coupling parameters are enough small was removed.

Some recent results dealt with the case when there is only kinematic dissipation or stress tensor
dissipation. Zhu \cite{Zhu2018} showed the existence of small global smooth solutions of the 3D Oldroyd-B model with $\eta  = 0$ and $a  = 0$. Inspired by the works of Zhu \cite{Zhu2018} and Danchin \cite{Danchin2001}, Chen and Hao \cite{ChenHao2019} extended such small data global
solution obtained in Sobolev spaces to the critical Besov spaces.
Zhai \cite{Zhai2021} established the existence of global solutions for highly oscillating initial data. For the case $\nu=0$,
when the damping mechanism and the Laplacian dissipation for $\tau$ exist, Elgindi and Rousset \cite{ElgindiRousset2015} established a global large solution in a certain sense by building a new quantity. Later, Liu and Elgindi extended
these results to the 3D in \cite{ElgindiLiu2015}. More well-posedness results can be found in \cite{HieberWenZi2019,HuangWangWenZi2022,ZiFangZhang2014}.

In this paper, we consider the following Oldroyd-B model  with
fractional stress tensor diffusion: 
	\begin{align}\label{eq2}
         \left\{\begin{array}{l}
			\partial_tu + u\cdot\nabla u+\nabla {\rm P} = {\rm div}~\tau,~~~~{\rm div}~u=0,\\[1ex]
			\partial_t\tau + u\cdot\nabla\tau+Q(\nabla u,\tau)+(-\Delta)^{\beta}\tau=D(u),\\[1ex]
			u|_{t=0}=u_0,~~\tau|_{t=0}=\tau_0, \\[1ex]
		\end{array}\right.
	\end{align}
where $\beta>0$ and fractional Laplacian operator $(-\Delta)^{\beta}$ is defined by
	$$
	(-\Delta)^{\beta}\tau=\mathscr{F}^{-1}(|\xi|^{2\beta}\widehat{\tau}),\quad\quad \widehat{\tau}(\xi,t)=\mathscr{F}\tau(\xi,t)=\int_{\mathbb{R}^d}e^{-ix\cdot\xi}\tau(x,t) dx.
	$$
For the model with fractional stress diffusion and no stress damping,
Constantin, Wu, Zhao and Zhu \cite{ConstantinWuZhaoZhu2021} proved global
well-posedness for sufficiently small initial data in $H^r(\mathbb R^d)$,
where $d=2,3$, $r>1+d/2$, and $1/2\leq\beta\leq1$. Substantial progress has also been made in critical Besov spaces. Wu and Zhao \cite{WuZhao2022} established small-data global
well-posedness in a hybrid homogeneous Besov framework adapted to the
different behavior of low and high frequencies. Chen, Liu, Qin and Ye \cite{ChenLiuQinYe2022} proved the analogous results in the $L^p$ framework. More recently, Liang, Luo and Yin \cite{liangluoyin} established small-data global well-posedness in
$H^s(\mathbb{R}^2)$ for $s \ge2\beta$ when $1/2 < \beta < 1$. These results provide motivation for studying Sobolev
regularity in three dimensions.

Now, we state our main result.
\begin{theo}\label{thm1}
Let $d=3$, $\beta\in(\frac12,1)$, and
$s\in[2\beta+\frac12,\frac52]$.
Assume that $(u_0,\tau_0)\in H^s\times H^s$,
$\nabla\cdot u_0=0$, and $\tau_0$ is symmetric.
Define
\[
E_0(t)\triangleq \|(u(t),\tau(t))\|_{H^\beta}^2 ~~~~\text{and}
~~~
E_1(t)\triangleq
\|\Lambda^\beta(u(t),\tau(t))\|_{H^{s-\beta}}^2.
\]
For every $M>0$, there exists $\epsilon=\epsilon(M)>0$
such that, if
\begin{equation}\label{cond1}
    E_0(0)\le M
    \qquad\text{and}\qquad
    E_1(0)\le\epsilon(M),
\end{equation}
then the following statements hold:

    (i) For $\beta\in(\frac12,1)$, there exists $T>0$
    such that \eqref{eq2} admits a unique local solution
    \[
    (u,\tau)\in C([0,T);H^s)\cap L_T^2([0,T);H^{s-\beta+1})\times C([0,T);H^s)\cap L_T^2([0,T);H^{s+\beta}).
    \]
    
 Moreover, the solution depends continuously on the initial data.

    (ii) If, in addition, $\beta\in(\frac12,\frac56]$, then the solution
    extends uniquely to a global solution
    \[
    (u,\tau)\in C([0,\infty);H^s)\times C([0,\infty);H^s)
    \]
\end{theo}
\begin{rema}
Unlike the classical local theory for the Euler equations in
$H^s$ with $s>\frac52$, our energy does not control
$\|\nabla u\|_{L^\infty}$. At $s=2\beta+\frac12$, the transport estimate takes the form
\[
\bigl|\langle\Lambda^s(u\cdot\nabla u),\Lambda^s u\rangle\bigr|
\lesssim
\|\Lambda^{2\beta+\frac{1}{2}} u\|_{L^2}
\|\Lambda^{\beta+\frac{3}{2}}u\|_{L^2}^2.
\]
This estimate uses the full dissipation and yields no additional
small factor from time integration..
Moreover, even when $s>2\beta+\frac12$, estimating
$\langle\Lambda^s Q(\nabla u,\tau),\Lambda^s\tau\rangle$
produces the contribution
\[
\|\tau\|_{L^\infty}
\|\Lambda^{s+1-\beta}u\|_{L^2}
\|\Lambda^{s+\beta}\tau\|_{L^2}.
\]
This also uses the full velocity and stress dissipation, with
no additional small factor from time integration. Thus, some smallness is
necessary to absorb these nonlinear contributions.
\end{rema}
\begin{rema}
    In the case $d=3$, Constantin, Wu, Zhao and Zhu \cite{ConstantinWuZhaoZhu2021} proved global
well-posedness for sufficiently small initial data in $H^s$ with $s>\frac{5}{2}$. Theorem \ref{thm1} relaxes the regularity index to $s\in [2\beta+\frac{1}{2},\frac{5}{2}]$.
\end{rema}
\begin{rema}
    We are concerned with the well-posedness in the sense of Hadamard, that is, the existence, uniqueness, and continuous dependence of the solution. By continuous dependence we mean that for any $(u_0,\tau_0)$ in some metric space $\mathbb{X}$, there exists a neighborhood $U$ of $(u_0,\tau_0)$ and $T=T(U)>0$ such that the solution map $(u_0,\tau_0)\to (u,\tau)$ is continuous from $U$ to $C([0,T];\mathbb{X})$.
\end{rema}
\begin{rema}
The restriction $\beta\le \frac56$ in Theorem~1.1(ii)
arises from the low-frequency derivative balance in the
estimates used in Section~4. Indeed, the velocity transport
term in the $\dot H^\beta$ energy satisfies
\[
\big|\langle \Lambda^\beta(u\cdot\nabla u),
                 \Lambda^\beta u\rangle\big|
\lesssim
\|\Lambda^{\frac52-2\beta}u\|_{L^2}
\|\Lambda^{2\beta}u\|_{L^2}^2.
\]
To control the first factor directly by
$\|\Lambda^\beta(u,\tau)\|_{H^{s-\beta}}$, we require
$\frac52-2\beta\ge\beta$, which is equivalent to
$\beta\le\frac56$.

The same restriction also appears in the closure of the lower-order energy in Section 4 where we obtain that
\[
|\langle Q(\nabla u,\tau),\tau\rangle|
\lesssim
\|\Lambda^\beta u\|_{L^2}
\|\Lambda^{\frac{5-2\beta}{4}}\tau\|_{L^2}^2.
\]
To control the last factor by the stress dissipation
$\|\Lambda^\beta\tau\|_{H^\beta}^2$, we require $\frac{5-2\beta}{4}\ge \beta$ which is equivalent to $\beta \le \frac{5}{6}$ as well. 
\end{rema}
\begin{coro}{\label{coro1.2}}{(Global solution with small data in critical case)}
   Let $d=3$, $\beta \in (\frac{1}{2},1)$ and $s=2\beta+\frac{1}{2}$.
   Assume that $(u_0,\tau_0)\in H^s\times H^s$,
$\nabla\cdot u_0=0$, and $\tau_0$ is symmetric. Then, there exists a small constant $\epsilon>0$ such that if 
    \begin{align*}
        \|(u_0,\tau_0)\|_{H^s}\le \epsilon,
    \end{align*}
    then the system (\ref{eq2}) has a unique global solution $(u,\tau)$ satisfying, for some constant $C>0$ and all $t>0$,
    \begin{align*}
        \|(u,\tau)\|_{H^s}\le C\epsilon.
    \end{align*}
    
\end{coro}
\begin{rema}
    The regularity threshold \(s=2\beta+\frac12\) is dictated by the balance between the velocity transport nonlinearity and the dissipation induced by the coupling. Indeed, the \(H^s\) energy yields velocity dissipation at order \(s+1-\beta\), while the transport term satisfies
\[
\big|\langle\Lambda^s(u\cdot\nabla u),\Lambda^s u\rangle\big|
\lesssim
\|\Lambda^{2\beta+\frac12}u\|_{L^2}
\|\Lambda^{s-\beta+1}u\|_{L^2}^{2}.
\]Thus, \(s=2\beta+\frac{1}{2}\) is the endpoint at which the coefficient of the dissipation can be controlled directly by the \(H^s\) energy and we use “critical” in this energy-estimate sense.
\end{rema}
\begin{rema}
    The corollary is a three-dimensional generalization of the critical global regularity result of Liang, Luo and Yin\cite{liangluoyin}, in which they consider the case $d=2, s=2\beta$. 
\end{rema}
\begin{rema}
J. Wu and J. Zhao \cite{WuZhao2022} proved if 
 $$\tau_0\in\dot{B}^{\frac{d}{2}+1-2\beta}_{2, 1}\cap \dot{B}^{\frac{d}{2}}_{2, 1}, u_0\in\dot{B}^{\frac{d}{2}+1-2\beta }_{2, 1}\cap \dot{B}^{\frac{d}{2}+2\beta-1}_{2, 1},$$
 then (\ref{eq2}) has a unique global solution with small initial data. Corollary \ref{coro1.2} implies the global well-posedness in critical Sobolev space $H^{2\beta+\frac{1}{2}}$ for $d=3$. Note that $ B^{2\beta+\frac{1}{2}}_{2, 1}\hookrightarrow H^{2\beta+\frac{1}{2}}$. 
\end{rema}
Theorem 1.1(i) asserts local existence, uniqueness
and continuous dependence in the $H^s$ topology. The proof is an application of the Bona-Smith method \cite{GuoLiYin2019}\cite{Masmoudi2007}.

The main analytical difficulty is that $H^s(\mathbb R^3)$ does not
embed into $W^{1,\infty}(\mathbb R^3)$ in the range considered here.
Therefore, the transport terms cannot be controlled solely by the
instantaneous $H^s$ norm through a Lipschitz estimate. Moreover, the
term $Q(\nabla u,\tau)$ contains a derivative of the undiffused velocity.
The analysis must instead use the regularization induced by the
coupling, with commutator and product estimates that distribute
derivatives according to the available dissipation.

The remainder of the paper is organized as follows. Section~2 collects
the notation and estimates from harmonic analysis. Section~3 presents
the proof of Theorem~1.1(i) , addressing local construction, uniqueness,
continuous dependence. Section~4 establishes the global regularity. In the last section, we give the proof of Corollary \ref{coro1.2}.

\section{Preliminaries}
In this section we collect some preliminary definitions and lemmas.
For more details we refer the readers to \cite{BCD2011,KatoPonce1988,Nirenberg1959,StefanovWuXuYe2025}.

We agree that $f\lm g$ represents $f\leq Cg$ with a constant $C>0$. The symbol $\widehat{f}=\mathscr{F}(f)$ stands for the Fourier transform of $f$. Denote $\mathscr{F}^{-1}(f)$ the inverse Fourier transform of $f$. We also write $\Lambda =(-\Delta)^{\frac{1}{2}}$ denoting the Zygmund operator.
	

Let $\chi:\mathbb{R}^d\to[0,1]$ be a radial, non-negative,
smooth and radially decreasing function which is supported in
$\mathcal{B}\triangleq\{\xi:|\xi|\leq\frac{4}{3}\}$
and $\chi\equiv1$ for $|\xi|\leq\frac{5}{4}$.
Let
$\varphi(\xi)=\chi(\frac{\xi}{2})-\chi(\xi)$.
Then $\varphi$ is supported in the ring
$\mathcal{C}\triangleq
\{\xi\in\mathbb{R}^d:\frac{3}{4}\leq|\xi|\leq\frac{8}{3}\}$.
For $u\in\mathcal{S}'$, $q\in\mathbb{Z}$,
we define the Littlewood-Paley operators:
$\dot{\Delta}_q u=\mathscr{F}^{-1}(\varphi(2^{-q}\cdot)\mathscr{F}u)$,
$\Delta_q u=\dot{\Delta}_q u$ for $q\geq0$,
$\Delta_q u=0$ for $q\leq-2$ and
$\Delta_{-1}u=\mathscr{F}^{-1}(\chi\mathcal{F}u)$,
and
$S_q u=\mathscr{F}^{-1}(\chi(2^{-q}\cdot)\mathcal{F}u)$.

We define the standard vector-valued Besov spaces
$B_{p,r}^s$ and $\dot{B}_{p,r}^s$
of the functions $u:\mathbb{R}^d\to\mathbb{R}^d$
with finite norms which are defined by
\begin{equation}
\|u\|_{B_{p,r}^s}
\triangleq
\left\|
\left(2^{js}\|\Delta_j u\|_{L^p}\right)_{j\in\mathbb{Z}}
\right\|_{\ell^r},
\end{equation}
\begin{equation}
\|u\|_{\dot{B}_{p,r}^s}
\triangleq
\left\|
\left(2^{js}\|\dot{\Delta}_j u\|_{L^p}\right)_{j\in\mathbb{Z}}
\right\|_{\ell^r}.
\end{equation}
Next we recall homogeneous Bony's decomposition
from \cite{BCD2011}.
Let $u$ and $v$ be tempered distributions in $\mathscr{S}'$, then the homogeneous $\rm paraproduct$ of $v$ by $u$ is defined as follows:
		\begin{align*}
			\dot{T}_u v = \sum_{j} \dot{S}_{j-1}u\dot{\Delta}_j v,
		\end{align*}
		and the homogeneous $\rm remainder$ of $u$ and $v$ is defined by
		\begin{align*}
			\dot{R}(u,v) = \sum_{|k-j|\leq1} \dot{\Delta}_{k} u \dot{\Delta}_{j} v.
		\end{align*}
		At least formally, we obtain the so-called Bony's decomposition:
		\begin{align*}
        uv=\dot{T}_u v + \dot{T}_v u + \dot{R}(u,v).
		\end{align*}

	\begin{lemm}\cite{BCD2011}\label{lemma2}	    	
		There exists a constant $C$ such that for any real number $s$ and any $(p,r)$ in $[1,\infty]^2$, we have, for any $(u,v)$ in $L^{\infty} \times \dot{B}^s_{p,r}$,
		\begin{align*}
			\|\dot{T}_u v\|_{\dot{B}^s_{p,r}} \leq C^{1+|s|} \|u\|_{L^{\infty}}\|v\|_{\dot{B}^s_{p,r}}.
		\end{align*}
		Moreover, for any $(s,t)$ in $\mathbb{R} \times (-\infty,0)$ and any $(p,r_1,r_2)$ in $[1,\infty]^3$, we have, for any $(u,v)\in \dot{B}^t_{\infty,r_1}\times\dot{B}^s_{p,r_2}$,
		\begin{align*}
			\|\dot{T}_u v\|_{\dot{B}^{s+t}_{p,r}} \leq \frac{C^{1+|s+t|}}{-t} \|u\|_{\dot{B}^{t}_{\infty,r_1}}\|v\|_{\dot{B}^s_{p,r_2}},
		\end{align*}
		with $\frac{1}{r}=\min\{1,\frac{1}{r_1}+\frac{1}{r_2}\}$.
	\end{lemm}
	\begin{lemm}\cite{BCD2011}\label{lemma3}
		A constant $C$ exists which satisfies the following inequalities. Let $(s_1,s_2)\in\mathbb{R}^2$ and $(p_1,p_2,r_1,r_2) \in [1,\infty]^4$. Assume that $$\frac{1}{p}=\frac{1}{p_1}+\frac{1}{p_2}\leq1~~~~~\text{and}~~~~~\frac{1}{r}=\frac{1}{r_1}+\frac{1}{r_2}\leq1.$$
		If $s_1+s_2>0$, for any $(u,v)\in \dot{B}^{s_1}_{p_1,r_1}\times \dot{B}^{s_2}_{p_2,r_2}$, then we have
		\begin{align*}
			\|\dot{R}(u,v)\|_{\dot{B}^{s_1+s_2}_{p,r}} \leq \frac{C^{1+|s_1+s_2|}}{s_1+s_2} \|u\|_{\dot{B}^{s_1}_{p_1,r_1}}\|v\|_{\dot{B}^{s_2}_{p_2,r_2}}.
		\end{align*}
	\end{lemm}
We introduce the Gagliardo-Nirenberg inequality for
$d = 3$.
\begin{lemm}\label{lemma4}\cite{Nirenberg1959}
For $d=3$, $p\in[2,+\infty)$ and $0\leq s,s_{1}\leq s_{2}$, there holds
\begin{align*}
\|\Lambda^{s}f\|_{L^{p}}\leq C\|\Lambda^{s_{1}}f\|_{L^{2}}^{1- \theta}\|\Lambda^{s_{2}}f\|_{L^{2}}^{\theta},
\end{align*}
where $0\leq\theta\leq 1$ and satisfies  
$$s+\frac{3}{2}-\frac{3}{p}=(1-\theta)s_{1}+\theta s_{2}.$$ 
Note that we also require that $0<\theta<1, 0\leq s_{1}\leq s$, when $ p=\infty$.
\end{lemm}

The following commutator estimate and product estimate are useful to control convective term.
\begin{lemm}\label{lemma5}\cite{KatoPonce1988}
Assume that $s>0$, $p, p_1, p_4 \in (1,\infty)$ and $\frac{1}{p}=\frac{1}{p_1}+\frac{1}{p_2}=\frac{1}{p_3}+\frac{1}{p_4}$, then we obtain
\begin{align*}
\|[\Lambda^{s},f]g\|_{L^{p}}\leq C\left(\|\Lambda^{s}f\|_{L^{p_{1}}}\|g\|_{L^{p_{2} }}+\|\nabla f\|_{L^{p_{3}}}\|\Lambda^{s-1}g\|_{L^{p_{4}}}\right).
\end{align*}
\end{lemm}
\begin{lemm}\label{lemma6}\cite{KatoPonce1988}
Assume that $s>0$, $p, p_2, p_4 \in (1,\infty)$ and $\frac{1}{p}=\frac{1}{p_1}+\frac{1}{p_2}=\frac{1}{p_3}+\frac{1}{p_4}$, then we obtain
\begin{align*}
\|\Lambda^{s}\left(fg\right)\|_{L^{p}} \leq C\left(\|f \|_{L^{p_{1}}}\|\Lambda^{s}g\|_{L^{p_{2}}} + \|g\|_{L^{p_ {3}}}\|\Lambda^{s}f\|_{L^{p_{4}}}\right).
\end{align*}
\end{lemm}
\begin{lemm}\label{lemma7}\cite{StefanovWuXuYe2025}
Assume that $0<s<1$, $p, r_1, r_2 \in (1,\infty)$ and $\frac{1}{p}=\frac{1}{r_1}+\frac{1}{r_2}$, then we obtain
\begin{align*}
\|[\Lambda^{s},f\cdot \nabla]g\|_{L^{p}}\leq C\|\nabla f\|_{L^{r_{1}}}\|\Lambda^{s}g\|_{L^{r_{2}}}.
\end{align*}
\end{lemm}
\section{Proof of Theorem 1.1 (i)}
	\par
    For convenience, we adopt the following notation  
    \begin{align*}
        E_0^{n}&\triangleq\|(u^n,\tau^n)\|_{H^{\beta}}^2+ 2k\langle u^n,\mathrm{div}\tau^n\rangle , \\
        D_0^n&\triangleq\frac{k}{2}\|\nabla u^n\|_{L^2}^2+\|\Lambda^{\beta }\tau^n\|^2_{H^{\beta}} ,\\ 
        E_1^n&\triangleq\|\Lambda^{\beta}(u^n,\tau^n)\|^2_{H^{s-\beta}}
        +2k\langle \Lambda^{2\beta-1}u^n,\Lambda^{2\beta-1}\mathrm{div}\tau^n \rangle+2k\langle \Lambda^{s-\beta}u^n,\Lambda^{s-\beta}\mathrm{div}\tau^n \rangle ,\\ 
        D_1^n&\triangleq \frac{k}{2}\|\Lambda^{2\beta} u^n\|^2_{H^{s-3\beta+1}}+\|\Lambda^{2\beta}\tau^n\|^2_{H^{s-\beta}}.
    \end{align*}
    
\begin{proof}

We first prove the local well-posedness. 
\\
Let $(u^{0},\tau^{0})=(0,0)$. We define a sequence $(u^{n},\tau^{n})_{n\in \mathbb{N}}$ by solving the following
linear equations:
\begin{align}\label{eq3}
\left\{\begin{array}{l}
\partial_tu^{n+1}+u^{n}\cdot\nabla u^{n+1}+\nabla P^{n+1} ={\rm div}~\tau^{n+1},~~~~{\rm div}~u^{n+1}=0,\\[1ex]
\partial_t\tau^{n+1}+u^{n}\cdot\nabla\tau^{n+1}+Q(\nabla u^{n},\tau^{n+1})+(-\Delta)^{\beta}\tau^{n+1}=D(u^{n+1}),\\[1ex]
u^{n+1}|_{t=0}=u_0,~~\tau^{n+1}|_{t=0}=\tau_0. \\[1ex]
\end{array}\right.
\end{align}

\textbf{Step 1: Uniform bounds}\\
Under the initial assumption of Theorem \ref{thm1}, for any constant $M>0$, there exists a constant $\epsilon(M)>0$ and $T>0$ such that for any $n\in \mathbb{N}$, we have 
\begin{align*}
   E_0^n(0) \le 2M \ , \ \ E_1^n(0) \le 2\epsilon. 
\end{align*}
Assume further that 
\begin{align}
    \label{3-2}&\sup_{t\in [0,T]} E_0^n(t)+\int_0^T D_0^n(t)dt\le 2M,\\
    \label{3-3}&\sup_{t\in [0,T]} E_1^n(t)+\int_0^T D_1^n(t)dt\le 2\epsilon.
\end{align}
Then we will prove the following estimates
\begin{align}
    \label{3-4}&\sup_{t\in [0,T]} E_0^{n+1}(t)+\int_0^T D_0^{n+1}(t)dt\le 2M,\\
    \label{3-5}&\sup_{t\in [0,T]} E_1^{n+1}(t)+\int_0^T D_1^{n+1}(t)dt\le 2\epsilon.
\end{align}

\textbf{(i)The energy estimate of $E_0^{n+1}$}
\\
Dotting (\ref{eq3}) by $(u^{n+1},\tau^{n+1})$ yields
\begin{align}{\label{3-6}}
&\quad\frac{1}{2}\frac{d}{dt}\|( u^{n+1},\tau^{n+1})\|_{L^{2}}^{2}+\|\Lambda^{\beta}\tau^{n+1}\|_{L^{2}}^{2}=-\langle Q(\nabla u^n,\tau^{n+1}),\tau^{n+1}\rangle.
\end{align}
Applying $\mathbb{P}$ to the first equation of (\ref{eq3}) and $\nabla\cdot$ to the second equation of (\ref{eq3}) separately, we have 
\begin{align} 
\label{3-7}&\partial_tu^{n+1}+\mathbb{P}(u^{n}\cdot\nabla u^{n+1}) =\mathbb{P}{\rm{div}}~\tau^{n+1},\\
\label{3-8}&\partial_t{\rm div}\tau^{n+1}+{\rm div}(u^{n}\cdot\nabla\tau^{n+1})+{\rm div}Q(\nabla u^{n},\tau^{n+1})+{\rm div}(-\Delta)^{\beta}\tau^{n+1}=\frac{1}{2}\Delta u^{n+1},
\end{align}
where $\mathbb{P}$ denotes the Leray projection.
Taking the $L^2$ inner product of (\ref{3-7}) with $k\rm{div}\tau^{n+1}$ and the $L^2$ inner product of (\ref{3-8}) with $ku^{n+1}$ yields
\begin{align}{\label{3-9}}
    &\frac{d}{dt}k\langle u^{n+1},\rm{div}\tau^{n+1}\rangle+\frac{k}{2}\|\nabla u^{n+1}\|_{L^2}^2 \\ \nonumber
    =&-k\langle\mathbb{P}(u^n \cdot \nabla u^{n+1}),\rm{div}\tau^{n+1}\rangle
      -k\langle \rm{div}(u^n\cdot \nabla \tau^{n+1}),u^{n+1}\rangle
      -k\langle\rm{div}Q(\nabla u^n,\tau^{n+1}),u^{n+1}\rangle\\ \nonumber
 &+k\langle\mathbb{P}\rm{div}\tau^{n+1},\rm{div}\tau^{n+1} \rangle
      -k\langle\rm{div}(-\Delta)^{\beta}\tau^{n+1},u^{n+1} \rangle .\nonumber
\end{align}
Applying $\Lambda^{\beta}$ to (\ref{eq3}) and then dotting by $(\Lambda^{\beta}u^{n+1},\Lambda^{\beta}\tau^{n+1})$ lead to
\begin{align}{\label{3-10}}
&\quad\frac{1}{2}\frac{d}{dt}\|\Lambda^{\beta}( u^{n+1},\tau^{n+1})\|_{L^{2}}^{2}+\|\Lambda^{2\beta}\tau^{n+1}\|_{L^{2}}^{2}\\\nonumber
&=-\langle \Lambda^{\beta}(u^n \cdot \nabla u^{n+1}),\Lambda^{\beta}u^{n+1}\rangle
-\langle\Lambda^{\beta} (u^n \cdot \nabla \tau^{n+1}),\Lambda^{\beta}\tau^{n+1}\rangle
-\langle \Lambda^{\beta}Q(\nabla u^n,\tau^{n+1}),\Lambda^{\beta}\tau^{n+1}\rangle.
\end{align}
Combining (\ref{3-6}),(\ref{3-9}) and (\ref{3-10}), we have
\begin{align}{\label{3-11}}
    &\frac{1}{2}\frac{d}{dt}E_0^{n+1}+D_0^{n+1} \\ \nonumber
    =& 
-\langle Q(\nabla u^n,\tau^{n+1}),\tau^{n+1}\rangle
  -k\langle\mathbb{P}(u^n \cdot \nabla u^{n+1}),\rm{div}\tau^{n+1}\rangle
      -k\langle \rm{div}(u^n\cdot \nabla \tau^{n+1}),u^{n+1}\rangle\\\nonumber
      &-k\langle\rm{div}Q(\nabla u^n,\tau^{n+1}),u^{n+1}\rangle      
  -\langle \Lambda^{\beta}(u^n \cdot \nabla u^{n+1}),\Lambda^{\beta}u^{n+1}\rangle
-\langle\Lambda^{\beta} (u^n \cdot \nabla \tau^{n+1}),\Lambda^{\beta}\tau^{n+1}\rangle\\\nonumber
&-\langle \Lambda^{\beta}Q(\nabla u^n,\tau^{n+1}),\Lambda^{\beta}\tau^{n+1}\rangle
 -k\langle\rm{div}(-\Delta)^{\beta}\tau^{n+1},u^{n+1} \rangle+k\langle\mathbb{P}\rm{div}\tau^{n+1},\rm{div}\tau^{n+1} \rangle \\\nonumber 
\triangleq&\sum_{i=1}^{9} J_i.
\end{align}
Transferring the derivatives, together with Lemma \ref{lemma4}, yields
\begin{align*}
     \sum_{i=1}^4 |J_i|   
  &\lesssim \|\nabla u^n\|_{L^2}\|\tau^{n+1}\|_{L^3}\|\tau^{n+1}\|_{L^6}+\|u^n\|_{L^{\infty}}\|\nabla u^{n+1}\|_{L^2}\|\nabla\tau^{n+1}\|_{L^2}\\[-0.9em]
  &\quad+\|\nabla u^n\|_{L^3}\|\tau^{n+1}\|_{L^6}\|\nabla u^{n+1}\|_{L^2}\\
   & \lesssim(E_1^n)^{\frac{1}{2}}\left(E_0^{n+1}+D_0^{n+1}\right).
\end{align*}
From Lemma \ref{lemma7}, we deduce that
\begin{align*}
    |J_5|\lesssim \|\nabla u^n\|_{L^{\frac{3}{2-2\beta}}}\|\Lambda^{\beta}u^{n+1}\|_{L^{\frac{6}{1+2\beta}}}^2
    \lesssim\|\Lambda^{2\beta+\frac{1}{2}}u^n\|_{L^2}\|\nabla u^{n+1}\|_{L^2}^2
    \lesssim (E_1^n)^{\frac{1}{2}} D_0^{n+1}.
\end{align*}
Using Lemma \ref{lemma4} and $\frac{1}{2}<\beta<1$, we obtain that
\begin{align*}
     |J_6|+|J_7|   
  &\lesssim \|u^n\|_{L^{\infty}}\|\nabla\tau^{n+1}\|_{L^2}\|\Lambda^{2\beta} \tau^{n+1}\|_{L^2}
   +\|\nabla u^n\|_{L^3}\|\tau^{n+1}\|_{L^6}\|\Lambda^{2\beta} \tau^{n+1}\|_{L^2}\\
   & \lesssim(E_1^n)^{\frac{1}{2}}D_0^{n+1}.
     \end{align*}
Finally, according to Hölder's and Young's inequalities, we have
\begin{align*}
    |J_8|&\le k\|\Lambda^{2\beta}\tau^{n+1}\|_{L^2}\|\nabla u^{n+1}\|_{L^2}\le Ck^2\|\nabla u^{n+1}\|_{L^2}^2+\frac{1}{100}\|\Lambda^{2\beta}\tau^{n+1}\|_{L^2}^2,\\
    |J_9|&\le Ck\|\nabla\tau^{n+1}\|_{L^2}^2 \ .
\end{align*}
Inserting the estimates above for $J_1$ through $J_{9}$ in (\ref{3-11}) and choosing $k$ sufficiently small, together with (\ref{3-2}) and (\ref{3-3}), we arrive at the following estimate
\begin{align}
     \frac{1}{2}\frac{d}{dt}E_0^{n+1}+(1-C\epsilon^{\frac{1}{2}})D_0^{n+1}
     \le C\epsilon^{\frac{1}{2}}E_0^{n+1}.
\end{align}
Define 
\begin{align*}
    A_0^n(t)\triangleq E_0^{n}(t)+\int_0^t D_0^{n}(t')dt'.
\end{align*}
It follows that 
\begin{align*}
     &\frac{1}{2}E_0^{n+1}(t)+(1-C\epsilon^{\frac{1}{2}}) \int_0^t D_0^{n+1}(t') dt'
     \le\frac{1}{2}E_0^{n+1}(0)+ CT\epsilon^{\frac{1}{2}}\sup_{t\in [0,T]}A_0^{n+1}(t).
\end{align*}
Taking $T^2\le \epsilon \le \min\{1,\frac{1}{16C^2}\}$, for any $t\in [0,T]$, one can get
\begin{align}
    \sup_{t\in [0,T]}A_0^{n+1}(t)\le 2M,
\end{align}
where we used $E_0^{n+1}(0)\le M$ which comes from (\ref{cond1}).

\textbf{(ii)The energy estimate of $E_1^{n+1}$}
\\
Applying $\Lambda^{s}$ to (\ref{eq3}) and then dotting by $(\Lambda^{s}u^{n+1},\Lambda^{s}\tau^{n+1})$ lead to
\begin{align}{\label{3-14}}
&\quad\frac{1}{2}\frac{d}{dt}\|\Lambda^{s}( u^{n+1},\tau^{n+1})\|_{L^{2}}^{2}+\|\Lambda^{s+\beta}\tau^{n+1}\|_{L^{2}}^{2}\\\nonumber
&=-\langle \Lambda^{s}(u^n \cdot \nabla u^{n+1}),\Lambda^{s}u^{n+1}\rangle
-\langle\Lambda^{s} (u^n \cdot \nabla \tau^{n+1}),\Lambda^{s}\tau^{n+1}\rangle
-\langle \Lambda^{s}Q(\nabla u^n,\tau^{n+1}),\Lambda^{s}\tau^{n+1}\rangle\\ \nonumber
&\triangleq\sum_{i=1}^3 J'_i.
\end{align}
Applying $\Lambda^{2\beta-1}$ to (\ref{3-7})-(\ref{3-8}) and taking $L^2$ inner product with $(k \Lambda^{2\beta-1}\rm{div} \tau^{n+1},k \Lambda^{2\beta-1}u^{n+1})$ yield 
\begin{align}{\label{3-15}}
  &\frac{d}{dt}k\langle \Lambda^{2\beta-1}u^{n+1},\Lambda^{2\beta-1}\mathrm{div}\tau^{n+1}\rangle_{L^{2}}+\frac{k}{2}\|\Lambda^{2\beta} u^{n+1}\|_{L^{2}}^{2} 
  \\\nonumber
  =&-k\langle\Lambda^{2\beta-1}(u^{n}\cdot\nabla u^{n+1}),\Lambda^{2\beta-1}\mathbb{P}\mathrm{div}\tau^{n+1}\rangle
  -k\langle\Lambda^{2\beta-1}\mathrm{div}(u^{n}\cdot\nabla\tau^{n+1}),\Lambda^{2\beta-1}u^{n+1}\rangle\\\nonumber
 &-k\langle\Lambda^{2\beta-1}\mathrm{div} Q(\nabla u^{n},\tau^{n+1}),\Lambda^{2\beta-1}u^{n+1}\rangle
  -k\langle\Lambda^{2\beta-1}(-\Delta)^\beta\mathrm{div}\tau^{n+1},\Lambda^{2\beta-1}u^{n+1}\rangle\\ \nonumber
  &+k\langle\mathbb{P}\Lambda^{2\beta-1}\mathrm{div}\tau^{n+1}, \Lambda^{2\beta-1}\mathrm{div}\tau^{n+1}\rangle\\ \nonumber
  \triangleq&\sum_{i=4}^8 J'_i.
\end{align}
Along the same line,
\begin{align}{\label{3-16}}
  &\frac{d}{dt}k\langle \Lambda^{s-\beta}u^{n+1},\Lambda^{s-\beta}\mathrm{div}\tau^{n+1}\rangle_{L^{2}}+\frac{k}{2}\|\Lambda^{s-\beta}\nabla u^{n+1}\|_{L^{2}}^{2} 
  \\\nonumber
  =&-k\langle\Lambda^{s-\beta}(u^{n}\cdot\nabla u^{n+1}),\Lambda^{s-\beta}\mathbb{P}\mathrm{div}\tau^{n+1}\rangle
  -k\langle\Lambda^{s-\beta}\mathrm{div}(u^{n}\cdot\nabla\tau^{n+1}),\Lambda^{s-\beta}u^{n+1}\rangle\\\nonumber
 &-k\langle\Lambda^{s-\beta}\mathrm{div} Q(\nabla u^{n},\tau^{n+1}),\Lambda^{s-\beta}u^{n+1}\rangle
  -k\langle\Lambda^{s-\beta}(-\Delta)^\beta\mathrm{div}\tau^{n+1},\Lambda^{s-\beta}u^{n+1}\rangle\\ \nonumber
  &+k\langle\mathbb{P}\Lambda^{s-\beta}\mathrm{div}\tau^{n+1}, \Lambda^{s-\beta}\mathrm{div}\tau^{n+1}\rangle\\ \nonumber
  \triangleq&\sum_{i=9}^{13} J'_i.
\end{align}
Combining (\ref{3-10}),(\ref{3-14}),(\ref{3-15}) and (\ref{3-16}), we have
\begin{align}{\label{3-17}}
    \frac{1}{2}\frac{d}{dt}E_1^{n+1}+D_1^{n+1} 
    = \sum_{i=5}^7 J_i +\sum_{i=1}^{13} J'_i.
\end{align}
It can be deduced directly that
\begin{align*}
     \sum_{i=5}^7 |J_i|   
 \lesssim (E_1^n)^{\frac{1}{2}}E_1^{n+1}.
\end{align*} 
By virtue of $\nabla\cdot u^n=0$ and Lemma \ref{lemma4}, we obtain
\begin{align*}
    |J'_1|&\lesssim\|\Lambda^s u^n\|_{L^{2}} \|\nabla u^{n+1}\|_{L^{\frac{3}{1-\beta}}}\|\Lambda^s u^{n+1}\|_{L^{\frac{6}{2\beta+1}}} + \|\nabla u^n\|_{L^{\frac{3}{2-2\beta}}} \|\Lambda^s u^{n+1}\|_{L^{\frac{6}{2\beta+1}}}^2\\ \nonumber
    &\lesssim\|\Lambda^{s} u^n\|_{L^{2}} \|\Lambda^{\beta+\frac{3}{2}} u^{n+1}\|_{L^{2}}\|\Lambda^{s-\beta+1} u^{n+1}\|_{L^{2}} + \|\Lambda^{2\beta+\frac{1}{2}} u^n\|_{L^{2}} \|\Lambda^{s-\beta+1} u^{n+1}\|_{L^2}^2\\ \nonumber
    &\lesssim(E_1^n)^{\frac{1}{2}} D_1^{n+1}.
\end{align*}
For $J'_2$, replacing $u^{n+1}$ by $\tau^{n+1}$ immediately yields the same estimate as for $J'_1$.
Using $s>\frac{3}{2}$ and $\| \tau^{n+1}\|_{L^{\infty}}\le \epsilon\|\tau^{n+1}\|_{L^2}+C\|\Lambda^s \tau^{n+1}\|_{L^2}$, we obtain that
 \begin{align*}
    |J'_3|&\lesssim\left(\|\Lambda^{s-\beta+1} u^n\|_{L^2} \|\tau^{n+1}\|_{L^{\infty}}+ \| \nabla u^n\|_{L^3} \|\Lambda^{s-\beta} \tau^{n+1}\|_{L^6}\right)\|\Lambda^{s+\beta} \tau^{n+1}\|_{L^2}\\ 
    &\le D_1^n( \epsilon E_0^{n+1}+C E_1^{n+1})+ C E_1^nD_1^{n+1}+\frac{1}{100}D_1^{n+1}.
\end{align*}
In view of Hölder's inequality and Lemma \ref{lemma4}, we infer that
\begin{align*}
    \sum_{i=4}^6 |J'_i| &\lesssim \|\Lambda^{2\beta-1}u^n\|_{L^6}\|\nabla u^{n+1}\|_{L^3}\|\Lambda^{2\beta}\tau^{n+1}\|_{L^2}
+\|u^n\|_{L^{\infty}}\|\Lambda^{2\beta}u^{n+1}\|_{L^2}\|\Lambda^{2\beta}\tau^{n+1}\|_{L^2}\\[-0.9em]
&\quad + \|\Lambda^{2\beta-1}u^n\|_{L^6}\|\nabla \tau^{n+1}\|_{L^3}\|\Lambda^{2\beta}u^{n+1}\|_{L^2}
+\|u^n\|_{L^{\infty}}\|\Lambda^{2\beta}\tau^{n+1}\|_{L^2}\|\Lambda^{2\beta}u^{n+1}\|_{L^2}\\
&\quad + \|\Lambda^{2\beta}u^n\|_{L^3}\|\tau^{n+1}\|_{L^6}\|\Lambda^{2\beta}u^{n+1}\|_{L^2}+\|\nabla u^n\|_{L^3}\|\Lambda^{2\beta-1}\tau^{n+1}\|_{L^6}\|\Lambda^{2\beta}u^{n+1}\|_{L^2}\\
&\lesssim(E_1^n)^{\frac{1}{2}} E_1^{n+1}.
\end{align*}
Young's inequality implies 
\begin{align*}
    |J'_7|&\le \frac{1}{100}\|\Lambda^{4\beta-1}\tau^{n+1}\|_{L^2}^2
    +Ck^2\|\Lambda^{2\beta} u^{n+1}\|_{L^2}^2,\\
    |J'_8| &\le Ck \|\Lambda^{2\beta}\tau^{n+1}\|_{L^2}^2.
\end{align*}
Thanks to $\frac{1}{2}<\beta<1$, we get
\begin{align*}
    |J'_{9}|& \lesssim \left(\|\Lambda^{s-\beta}u^n\|_{L^{\frac{6}{2\beta+1}}}\|\nabla u^{n+1}\|_{L^{\frac{3}{1-\beta}}}+\|u^n\|_{L^{\infty}}\|\Lambda^{s-\beta+1}u^{n+1}\|_{L^2}\right)\|\Lambda^{s-\beta+1}\tau^{n+1}\|_{L^2}\\
   &\lesssim\left(\|\Lambda^{s-2\beta+1} u^n\|_{L^2}\|\Lambda^{\beta+\frac{3}{2}} u^{n+1}\|_{L^2}+\|u^n\|_{L^{\infty}}\|\Lambda^{s-\beta+1} u^{n+1}\|_{L^2}\right)\|\Lambda^{s-\beta+1}\tau^{n+1}\|_{L^2}\\
    &\lesssim (E_1^n)^{\frac{1}{2}} D_1^{n+1}.
\end{align*}  
Obviously, $|J'_{10}|$ admits the same estimate as $|J'_{9}|$, with $u^{n+1}$ and $\tau^{n+1}$ interchanged.\\
Just like $J'_{3}$, we have
\begin{align*}
|J'_{11}|&\lesssim \left(\|\Lambda^{s-\beta+1} u^n\|_{L^2} \|\tau^{n+1}\|_{L^{\infty}}+ \| \nabla u^n\|_{L^3} \|\Lambda^{s-\beta} \tau^{n+1}\|_{L^6}\right)\|\Lambda^{s-\beta+1} u^{n+1}\|_{L^2} \\    
    &\le D_1^n( \epsilon E_0^{n+1}+C E_1^{n+1})+ C E_1^nD_1^{n+1}+\frac{1}{100}D_1^{n+1}.
\end{align*}
Moreover, one can derive directly that
\begin{align*}
    |J'_{12}|\le \frac{1}{100}\|\Lambda^{s+\beta}\tau^{n+1}\|_{L^2}^2+Ck^2\|\Lambda^{s-\beta+1}u^{n+1}\|_{L^2}^2,
\end{align*}
and
\begin{align*}
    |J'_{13}| \le Ck\|\Lambda^{s-\beta+1}\tau^{n+1}\|_{L^2}^2.
\end{align*}
Inserting the estimates above in (\ref{3-17}) and choosing $k$ sufficiently small, we arrive at the following estimate
\begin{align*}
    \frac{1}{2}\frac{d}{dt}E_1^{n+1}+\left(1-C(\epsilon+\epsilon^{\frac{1}{2}})\right)D_1^{n+1} 
    \le  C D_1^{n}(\epsilon E_0^{n+1}+C E_1^{n+1})+C \epsilon^{\frac{1}{2}}E_1^{n+1}.
\end{align*}
Define 
\begin{align*}
    A_1^n(t)\triangleq E_1^{n}(t)+\int_0^t D_1^{n}(t')dt'.
\end{align*}
It follows that
\begin{align}
    &\frac{1}{2}E_1^{n+1}(t)+\left(1-C(\epsilon+\epsilon^{\frac{1}{2}})\right)\int_0^t D_1^{n+1}(t')dt'\\ \nonumber
    \le&\frac{1}{2}E_1^{n+1}(0)+C\epsilon^2+C\left(\epsilon+T\epsilon^{\frac{1}{2}} \right)\sup_{t\in[0,T]}A_1^{n+1}(t),
\end{align}
where we used $\sup_{[0,T]}E_0^{n+1}$ is bounded.\\
Taking $T^{2}\le \epsilon \le \min \{1,\frac{1}{16C},\frac{1}{16C^2} \}$, for any $t\in [0,T]$, one can get
\begin{align}
    \sup_{t\in [0,T]}A_1^{n+1}(t)\le 2\epsilon.
\end{align}
\\

\textbf{Step 2: Convergence}
\\
Denote that $\bar{u}^{n}\triangleq u^{n+1}-u^n,\bar{\tau}^{n}\triangleq \tau^{n+1}-\tau^n$ and 
\begin{align*}
    \bar{E}^n&\triangleq\|(\bar{u}^n,\bar{\tau}^n)\|_{H^{s-\beta}}^2+2k\langle\bar{u}^n,\rm{div}\bar{\tau}^n\rangle_{H^{s-2\beta}},\\
    \bar{D}^n&\triangleq\frac{k}{2}\|\nabla \bar{u}^n\|_{H^{s-2\beta}}^2+\|\Lambda^{\beta}\bar{\tau}^n\|_{H^{s-\beta}}^2.
\end{align*}
According to (\ref{eq3}), we deduce that 
\begin{align}\label{eq4}
\left\{\begin{array}{l}
\partial_t\bar{u}^n+u^{n}\cdot\nabla \bar{u}^n+\nabla \bar{P}^n ={\rm div}~\bar{\tau}^n-\bar{u}^{n-1}\cdot\nabla u^n,~~~~{\rm div}~\bar{u}^n=0,\\[1ex]
\partial_t\bar{\tau}^n+u^{n}\cdot\nabla\bar{\tau}^n+Q(\nabla u^{n},\bar{\tau}^n)+Q(\nabla \bar{u}^{n-1},\tau^n)=D(\bar{u}^{n})-(-\Delta)^{\beta}\bar{\tau}^n-\bar{u}^{n-1}\cdot\nabla\tau^n,\\[1ex]
\bar{u}^n|_{t=0}=0,~~\bar{\tau}^n|_{t=0}=0. \\[1ex]
\end{array}\right.
\end{align}
Proceeding as in Step 1, we obtain that
\begin{align}
    &\quad\frac{1}{2}\frac{d}{dt}(\|\bar{u}^n\|_{H^{s-\beta}}^2+\|\bar{\tau}^n\|_{H^{s-\beta}}^2+2k\langle\bar{u}^n,\rm{div}\bar{\tau}^n\rangle_{H^{s-2\beta}})+\frac{k}{2}\|\nabla \bar{u}^n\|_{H^{s-2\beta}}^2+\|\Lambda^{\beta} \bar{\tau}^n\|_{H^{s-\beta}}^2\\ \nonumber
     &\triangleq\sum_{i=1}^{14} K_i
\end{align}
where
\begin{align*}
    &K_1=-\langle \bar{u}^{n-1}\cdot\nabla u^n,\bar{u}^n\rangle-\langle \bar{u}^{n-1}\cdot\nabla \tau^n,\bar{\tau}^n\rangle
-\langle Q(\nabla u^{n},\bar{\tau}^n),\bar{\tau}^n\rangle-\langle Q(\nabla \bar{u}^{n-1},\tau^n),\bar{\tau}^n\rangle\\
    &~~~~~~~~~~
    -k\langle \mathbb{P}(u^n\cdot\nabla \bar{u}^n),\mathrm{div}\bar{\tau}^n\rangle
    -k\langle \mathbb{P}(\bar{u}^{n-1}\cdot\nabla u^n),\mathrm{div}\bar{\tau}^n\rangle
    -k\langle \mathrm{div}(u^n\cdot\nabla\bar{\tau}^n),\bar{u}^n\rangle \\ 
    &~~~~~~~~~~-k\langle \mathrm{div}(\bar{u}^{n-1}\cdot \nabla\tau^n),\bar{u}^n\rangle
    -k\langle \mathrm{div}Q(\nabla u^{n},\bar{\tau}^n),\bar{u}^n\rangle
    -k\langle \mathrm{div}Q(\nabla \bar{u}^{n-1},\tau^n),\bar{u}^n\rangle,\\
    &K_{2}= -\langle\Lambda^{s-\beta}(u^n\cdot\nabla \bar{u}^n),\Lambda^{s-\beta} \bar{u}^n\rangle,~~~~~~~~~~
    K_3=-\langle \Lambda^{s-\beta}(\bar{u}^{n-1}\cdot\nabla u^n),\Lambda^{s-\beta}\bar{u}^n\rangle,\\
    &K_{4}=-\langle\Lambda^{s-\beta}(u^n\cdot\nabla \bar{\tau}^n),\Lambda^{s-\beta} \bar{\tau}^n\rangle,~~~~~~~~~~
    K_5=-\langle \Lambda^{s-\beta}(\bar{u}^{n-1}\cdot\nabla \tau^n),\Lambda^{s-\beta}\bar{\tau}^n\rangle,\\
    &K_{6}=-\langle \Lambda^{s-\beta}Q(\nabla u^{n},\bar{\tau}^n),\Lambda^{s-\beta}\bar{\tau}^n\rangle,~~~~~~~~
    K_{7}=-\langle \Lambda^{s-\beta}Q(\nabla \bar{u}^{n-1},\tau^n),\Lambda^{s-\beta}\bar{\tau}^n\rangle,\\
    &K_{8}= -k\langle \Lambda^{s-2\beta}\mathbb{P}(u^n\cdot\nabla \bar{u}^n),\Lambda^{s-2\beta}\mathrm{div}\bar{\tau}^n\rangle,~~~~~
    K_{9}=-k\langle \Lambda^{s-2\beta}\mathbb{P}(\bar{u}^{n-1}\cdot\nabla u^n),\Lambda^{s-2\beta}\mathrm{div}\bar{\tau}^n\rangle,\\
    &K_{10}=-k\langle \Lambda^{s-2\beta}\mathrm{div}(u^n\cdot\nabla\bar{\tau}^n),\Lambda^{s-2\beta}\bar{u}^n\rangle,~~~~~
     K_{11}=-k\langle \Lambda^{s-2\beta}\mathrm{div}(\bar{u}^{n-1}\cdot \nabla\tau^n),\Lambda^{s-2\beta}\bar{u}^n\rangle,\\
     &K_{12}=-k\langle \Lambda^{s-2\beta}\mathrm{div}Q(\nabla u^{n},\bar{\tau}^n),\Lambda^{s-2\beta}\bar{u}^n\rangle,~~~
    K_{13}=-k\langle \Lambda^{s-2\beta}\mathrm{div}Q(\nabla \bar{u}^{n-1},\tau^n),\Lambda^{s-2\beta}\bar{u}^n\rangle,\\
    &K_{14}=-k\langle \mathrm{div}\Lambda^{2\beta} \bar{\tau}^n,\bar{u}^n\rangle
    -k\langle \mathrm{div}\Lambda^{s} \bar{\tau}^n,\Lambda^{s-2\beta}\bar{u}^n\rangle
     +k\langle \mathbb{P}\mathrm{div}\bar{\tau}^n,\mathrm{div}\bar{\tau}^n\rangle_{H^{s-2\beta}}.
\end{align*}
We first note that  \[\|\nabla( u^n,\tau^n)\|_{L^2}\lesssim\epsilon^{\frac{1}{2}},~~
\|\nabla( u^n,\tau^n)\|_{L^3}\lesssim\epsilon^{\frac{1}{2}},~~\|(u^n,\tau^n)\|_{L^{\infty}}\lesssim\epsilon^{\frac{1}{2}}.\] 
Thus, owing to Hölder's and Sobolev's inequalities, it is not difficult to check that
\begin{align*}
     |K_1| \lesssim \epsilon^{\frac{1}{2}}\left(\bar{E}^{n-1}+\bar{E}^n\right).
\end{align*}
By $\frac{1}{2} < \beta<1$, we get 
\begin{align*}
    |K_{2}|&\lesssim\|\Lambda^{s-\beta}u^n\|_{L^{\frac{6}{3-2\beta}}}\|\nabla\bar{u}^n\|_{L^{3}}\|\Lambda^{s-\beta}\bar{u}^n\|_{L^{\frac{6}{2\beta+1}}}+\|\nabla u^n\|_{L^{\frac{3}{2-2\beta}}}\|\Lambda^{s-\beta}\bar{u}^n\|_{L^{\frac{6}{2\beta+1}}}^2\\
    &\lesssim \|\Lambda^{s}u^n\|_{L^{2}}\|\Lambda^{\frac{3}{2}}\bar{u}^n\|_{L^{2}}\|\Lambda^{s-2\beta+1}\bar{u}^n\|_{L^{2}}+\|\Lambda^{2\beta+\frac{1}{2}} u^n\|_{L^{2}}\|\Lambda^{s-2\beta+1}\bar{u}^n\|_{L^{2}}^2\\
    &\lesssim \epsilon^{\frac{1}{2}} \bar{D}^n.
\end{align*}
As $\beta+\frac{1}{2}\le s-\beta$, we infer that
\begin{align*}
	|K_{3}|&\lesssim\|\Lambda^{s-1}\bar{u}^{n-1}\|_{L^{\frac{6}{4\beta-1}}} \|\nabla u^n\|_{L^{\frac{3}{2-2\beta}}}  \|\Lambda^{s-2\beta+1}\bar{u}^n\|_{L^2}
	+\|\bar{u}^{n-1}\|_{L^{\frac{3}{1-\beta}}}  \|\Lambda^s u^n\|_{L^{\frac{6}{1+2\beta}}}  \|\Lambda^{s-2\beta+1}\bar{u}^n\|_{L^2}\\
	& \le C \left(\|\Lambda^{2\beta+\frac{1}{2}} u^n\|_{L^2}^2\|\Lambda^{s-2\beta+1}\bar{u}^{n-1}\|_{L^2}^2
	  +\|\Lambda^{\beta+\frac{1}{2}} \bar{u}^{n-1}\|_{L^{2}}^2\|\Lambda^{s-\beta+1}u^n\|_{L^2}^2\right)
     + \frac{k}{100} \|\Lambda^{s-2\beta+1}\bar{u}^n\|_{L^2}^2\\
	&\le C \epsilon^{\frac{1}{2}} \bar{E}^{n-1}+C D_1^n \bar{E}^{n-1}+\frac{k}{100} \|\Lambda^{s-2\beta+1}\bar{u}^n\|_{L^2}^2.
\end{align*}
A similar estimate holds for $K_{4}$ and $K_{5}$, that is
\begin{align*}
     |K_{4}|+|K_{5}|&\le C\epsilon^{\frac{1}{2}} \bar{D}^n + C\epsilon^{\frac{1}{2}} \bar{E}^{n-1}+C D_1^n \bar{E}^{n-1}+\frac{k}{100} \|\Lambda^{s-2\beta+1}\bar{u}^n\|_{L^2}^2.
\end{align*}
In view of $s >\frac{3}{2}$, we have 
\begin{align*}
	|K_{6}|&\lesssim \|\nabla u^n\|_{L^{3}} \|\Lambda^{s-2\beta}\bar{\tau}^n\|_{L^{6}} \|\Lambda^s \bar{\tau}^n\|_{L^2}+\|\Lambda^{s-2\beta+1}u^n\|_{L^2} \|\bar{\tau}^n\|_{L^{\infty}}\|\Lambda^s \bar{\tau}^n\|_{L^2}\\
	& \lesssim \|\Lambda^{\frac{3}{2}} u^n\|_{L^2}\|\Lambda^{s-2\beta+1}\bar{\tau}^n\|_{L^2}\|\Lambda^s \bar{\tau}^n\|_{L^2}+\|\Lambda^{s-2\beta+1}u^n\|_{L^2}
	(\|\bar{\tau}^n\|_{L^{2}}+\|\Lambda^s \bar{\tau}^n\|_{L^2})\|\Lambda^s \bar{\tau}^n\|_{L^2}\\
	& \lesssim \epsilon^{\frac{1}{2}}\left(\bar{E}^n+\bar{D}^n \right).
\end{align*}
Noting $\beta>\frac{1}{2}$ and $\frac{3}{2}\le s-2\beta+1$ yields
\begin{align*}
	|K_{7}|&\lesssim \|\nabla \bar{u}^{n-1}\|_{L^{3}} \|\Lambda^{s-2\beta}\tau^n\|_{L^{6}}\|\Lambda^s\bar{\tau}^n\|_{L^2}+\|\Lambda^{s-2\beta+1}\bar{u}^{n-1}\|_{L^2}\|\tau^n\|_{L^{\infty}}\|\Lambda^s\bar{\tau}^n\|_{L^2}\\
    &\lesssim \|\Lambda^{s-2\beta+1}\tau^n\|_{L^2}\|\Lambda^{\frac{3}{2}} \bar{u}^{n-1}\|_{L^2}\|\Lambda^s \bar{\tau}^n\|_{L^2} 
    +\|\tau^n\|_{L^{\infty}}\|\Lambda^{s-2\beta+1}\bar{u}^{n-1}\|_{L^2}\|\Lambda^s\bar{\tau}^n\|_{L^2}\\
    &\lesssim \epsilon^{\frac{1}{2}}\left(\bar{D}^{n-1}+\bar{D}^n\right).
\end{align*}
Similarly, we have
\begin{align*}
    |K_{8}|&\lesssim \|\Lambda^{s-2\beta} {u}^n\|_{L^{6}}\|\nabla\bar{u}^n\|_{L^{3}}\|\Lambda^{s-2\beta+1}\bar{\tau}^n\|_{L^2}+\|u^n\|_{L^{\infty}}\|\Lambda^{s-2\beta+1}\bar{u}^n\|_{L^2}\|\Lambda^{s-2\beta+1}\bar{\tau}^n\|_{L^2}\\
    &\lesssim \epsilon^{\frac{1}{2}} \bar{D}^n.
\end{align*}
Since $s-2\beta+1>\frac{5}{2}-2\beta$ implies $\|\Lambda^{\frac{5}{2}-2\beta}\bar{u}^{n-1}\|_{L^2}\le \|\bar{u}^{n-1}\|_{L^2}+\|\Lambda^{s-2\beta+1}\bar{u}^{n-1}\|_{L^2}$, then we get
\begin{align*}
    |K_{9}| &\lesssim \|\Lambda^{s-2\beta}\bar{u}^{n-1}\|_{L^6}\|\nabla u^n\|_{L^3}\|\Lambda^{s-2\beta+1}\bar{\tau}^n\|_{L^2}+\|\bar{u}^{n-1}\|_{L^{\frac{3}{2\beta-1}}} \|\Lambda^{s-2\beta+1} u^n\|_{L^{\frac{6}{5-4\beta}}} \|\Lambda^{s-2\beta+1}\bar{\tau}^n\|_{L^2}\\
    &\lesssim\|\Lambda^{\frac{3}{2}}u^n\|_{L^2}\|\Lambda^{s-2\beta+1}\bar{u}^{n-1}\|_{L^2}
    \|\Lambda^{s-2\beta+1}\bar{\tau}^n\|_{L^2}+ \|\Lambda^{\frac{5}{2}-2\beta}\bar{u}^{n-1}\|_{L^2} \|\Lambda^{s} u^n\|_{L^2}\|\Lambda^{s-2\beta+1}\bar{\tau}^n\|_{L^2}\\
    &\lesssim  \epsilon^{\frac{1}{2}} \left( \bar{E}^{n-1}+\bar{D}^{n-1}+\bar{D}^n \right) .
\end{align*}
$K_{10}$ and $K_{11}$ can be treated in the same way as $K_{8}$ and $K_{9}$, while $K_{12}$ and $K_{13}$ can be handled similar to $K_{6}$ and $K_{7}$. Consequently, we obtain
\begin{align*}
    \sum_{i=10}^{13} |K_i| \lesssim \epsilon^{\frac{1}{2}}\left( \bar{E}^{n-1}+\bar{D}^{n-1}+\bar{E}^{n}+\bar{D}^n\right).
\end{align*}
Hölder's and Young's inequalities immediately show that
\begin{align*}
    |K_{14}|\le Ck^2\|\nabla \bar{u}^n\|_{H^{s-2\beta}}^2 + \frac{1}{100} \|\Lambda^{\beta}\bar{\tau}^n\|_{H^{s-\beta}}^2.
\end{align*}
Therefore, by adding up the estimates for $K_1 -K_{14}$ and choosing $k$ sufficiently small, we conclude that
\begin{align}{\label{4-21}}
    &\quad \ \frac{1}{2}\frac{d}{dt} \bar{E}^n (t)+ \bar{D}^n (t) 
    \le C\epsilon^{\frac{1}{2}}\left(\bar{E}^n (t)+\bar{D}^n(t)+\bar{E}^{n-1} (t) +\bar{D}^{n-1}(t)\right)+C D_1^n \bar{E}^{n-1}.
\end{align}
For convenience, we give the following notation 
\begin{align*}
    \bar{A}^n(t) \triangleq \bar{E}^n(t)+\int_0^t \bar{D}^n(t')dt'.
\end{align*}
Integrating (\ref{4-21}) and choosing $T$ and $\epsilon$ sufficiently small, we obtain that
\begin{align*}
   \sup_{t\in[0,T]}\bar{A}^n(t) \le \frac{1}{2}\sup_{t\in[0,T]}\bar{A}^{n-1}.
\end{align*}
A standard compactness argument then yields the existence of a local solution to system (\ref{eq2}).  
  \\  

\textbf{Step 3: Uniqueness}
\\
Assume $(u_1,\tau_1)$ and $(u_2,\tau_2)$ are two solutions of system (\ref{eq2}) with the same initial data. Then $\delta u\triangleq u_1-u_2,\delta \tau\triangleq \tau_1-\tau_2$ satisfy the following equations
\begin{align}\label{eq5}
\left\{\begin{array}{l}
\partial_t\delta u+u_1\cdot\nabla \delta u+\delta u \cdot\nabla u_2 + \nabla \delta P =\rm{div}~\delta\tau,~~~~{\rm div}~\delta u=0,\\[1ex]
\partial_t\delta \tau+u_1\cdot \nabla \delta\tau+\delta u\cdot \nabla\tau_2+Q(\nabla\delta u,\tau_1)+Q(\nabla u_2,\delta\tau)+\Lambda^{2\beta}\delta\tau=D(\delta u),\\[1ex]
\delta u|_{t=0}=0,~~\delta \tau|_{t=0}=0. \\[1ex]
\end{array}\right.
\end{align}
Denote 
\begin{align*}
    E_{\delta}(t)&\triangleq \|\delta u\|_{H^{s-\beta}}^2 +\|\delta \tau\|_{H^{s-\beta}}^2\\
    D_{\delta}(t)&\triangleq \frac{k}{2}\|\nabla \delta u\|_{H^{s-2\beta}}^2+\|\Lambda^{\beta}\delta \tau\|_{H^{s-\beta}}^2.
\end{align*}
Repeating the procedure in Step 2 and using Gronwall's inequality, we derive that
\begin{align}{\label{4-23}}
   \sup_{t\in [0,T]}E_{\delta}(t)+\int_0^{T}D_{\delta}(t)dt \le C E_{\delta}(0)\exp{\int_0^T F(t)dt}
\end{align}
where $F(t)=C\left(\sum\limits_{i=1}^2 ~(\|(u_i,\tau_i)\|_{H^s}^2+\|(u_i,\tau_i)\|_{H^s}+\|\nabla u_i\|_{H^{s-\beta}}^2+\|\Lambda^{\beta}\tau_i\|_{H^s}^2)\right)$.\\
The uniqueness then follows since $E_{\delta}(0)=0$.
\\

\textbf{Step 4: Continuous dependence}
\\
Inspired by \cite{GuoLiYin2019} and \cite{Masmoudi2007}, we use the Bona-Smith method to prove the continuous dependence. For this purpose, we need higher-regularity estimate for the solution. More precisely, 
\begin{lemm}{(Propagation of higher regularity)}{\label{highregularity}}
    Let $\frac12<\beta<1$, $2\beta+\frac12\le s\le \frac52$ and
$r\ge s$. Let $(u,\tau)$ be the solution of system (\ref{eq2}) on $[0,T]$
constructed above. Assume that
\[
\sup_{t\in[0,T]}
\|\Lambda^\beta(u,\tau)(t)\|_{H^{s-\beta}}
\le\delta,
\qquad
\int_0^T
\left(
\|\nabla u\|_{H^{s-\beta}}^2
+\|\Lambda^\beta\tau\|_{H^s}^2
\right)dt<\infty.
\]
There exists $\delta>0$ such that,
if $(u_0,\tau_0)\in H^r\times H^r$,
then
\begin{align}
&\sup_{t\in[0,T]}\|(u,\tau)(t)\|_{H^r}^2
+\int_0^T
\left(
\|\nabla u\|_{H^{r-\beta}}^2
+\|\Lambda^\beta\tau\|_{H^r}^2
\right)dt \le
C\|(u_0,\tau_0)\|_{H^r}^2.
\label{eq:higher-regularity}
\end{align}
Here $C$ depends only on $\beta,s,r,T$. We stress that the constants C here and in what follows in this step are independent of N.
\end{lemm}
\begin{proof}
    For convenience, we adopt the following notation:
\begin{align*}
E_r
&\triangleq\|(u,\tau)\|_{L^2}^2
+\|\Lambda^r(u,\tau)\|_{L^2}^2
+2k\langle u,\operatorname{div}\tau\rangle
+2k\langle\Lambda^{r-\beta}u,
             \Lambda^{r-\beta}\operatorname{div}\tau\rangle,\\
 D_r
&\triangleq\frac{k}{2}
\left(
\|\nabla u\|_{L^2}^2
+\|\Lambda^{r-\beta+1}u\|_{L^2}^2
\right)
+\|\Lambda^\beta\tau\|_{L^2}^2
+\|\Lambda^{r+\beta}\tau\|_{L^2}^2.
\end{align*}
Indeed,
\begin{align}{\label{high1}}
    \sum_{a\in\{0,r-\beta\}}
\left|
\langle\Lambda^a u,\Lambda^a\operatorname{div}\tau\rangle
\right|
\lesssim
\sum_{a\in\{0,r-\beta\}}
\|\Lambda^{a+\frac12}(u,\tau)\|_{L^2}^2
\lesssim\|(u,\tau)\|_{H^r}^2.
\end{align}
By virtue of $\beta<\frac32<s$ and
$\beta<2\beta+\frac12\le s$, we obtain
\begin{align}{\label{high2}}
\|(u,\tau)\|_{L^\infty}
+\|\Lambda^{2\beta+\frac12}u\|_{L^2}
&\lesssim
\|\Lambda^\beta(u,\tau)\|_{H^{s-\beta}}
\lesssim\delta.
\end{align}
{\bf (i) The lower-order energy estimate.}
Taking the $L^2$ inner product of (\ref{eq2}) with $(u,\tau)$
and using the cancellation of the coupling terms yield
\[
\frac12\frac{d}{dt}\|(u,\tau)\|_{L^2}^2
+\|\Lambda^\beta\tau\|_{L^2}^2
=-\langle Q(\nabla u,\tau),\tau\rangle.
\]
Along the same line as (\ref{3-9}), we have
\begin{align*}
\frac{d}{dt}k\langle u,\operatorname{div}\tau\rangle
+\frac{k}{2}\|\nabla u\|_{L^2}^2
&=-k\langle\mathbb P(u\cdot\nabla u),
                    \operatorname{div}\tau\rangle
-k\langle\operatorname{div}(u\cdot\nabla\tau),u\rangle\\
&\quad
-k\langle\operatorname{div}Q(\nabla u,\tau),u\rangle
-k\langle\Lambda^{2\beta}\operatorname{div}\tau,u\rangle
+k\|\mathbb P\operatorname{div}\tau\|_{L^2}^2.
\end{align*}
Transferring the derivatives and using
(\ref{high2}), we obtain
\begin{align*}
&\quad \frac12\frac{d}{dt}\left(\|(u,\tau)\|_{L^2}^2+2k\langle u,\operatorname{div}\tau\rangle\right)
+\|\Lambda^\beta\tau\|_{L^2}^2+\frac{k}{2}\|\nabla u\|_{L^2}^2\\
&\lesssim
\|\tau\|_{L^\infty}\|\nabla u\|_{L^2}\|\tau\|_{L^2}
+k\|u\|_{L^\infty}
  \|\nabla u\|_{L^2}\|\nabla\tau\|_{L^2}\\
&\quad
+k\|\tau\|_{L^\infty}\|\nabla u\|_{L^2}^2
+k\|\Lambda^{2\beta}\tau\|_{L^2}\|\nabla u\|_{L^2}
+k\|\nabla\tau\|_{L^2}^2\\
&\lesssim
\delta E_r^{\frac{1}{2}}D_r^{\frac{1}{2}}
+k\delta D_r
+kD_r\\
&\le
\left(\frac1{100}+C\delta+Ck\right) D_r
+C\delta^2 E_r.
\end{align*}

\noindent
{\bf (ii) The higher-order energy estimate.}
Applying $\Lambda^r$ to (\ref{eq2}) and then taking the $L^2$
inner product with $(\Lambda^r u,\Lambda^r\tau)$ lead to
\begin{align*}
&\frac12\frac{d}{dt}
\|\Lambda^r(u,\tau)\|_{L^2}^2
+\|\Lambda^{r+\beta}\tau\|_{L^2}^2 \notag\\
=&-\langle\Lambda^r(u\cdot\nabla u),\Lambda^r u\rangle
-\langle\Lambda^r(u\cdot\nabla\tau),\Lambda^r\tau\rangle
-\langle\Lambda^rQ(\nabla u,\tau),\Lambda^r\tau\rangle
\notag\\
\triangleq&\mathcal{J}_1+\mathcal{J}_2+\mathcal{J}_3.
\label{eq:hr-high-energy}
\end{align*}
By virtue of $\operatorname{div}u=0$ and Lemma \ref{lemma5},
\begin{align*}
|\mathcal{J}_1|
&\lesssim
\|\nabla u\|_{L^{\frac3{2-2\beta}}}
\|\Lambda^r u\|_{L^{\frac6{1+2\beta}}}^2\\
&\lesssim
\|\Lambda^{2\beta+\frac12}u\|_{L^2}
\|\Lambda^{r-\beta+1}u\|_{L^2}^2
\lesssim\delta D_r.
\end{align*}
Lemma \ref{lemma6} gives
\begin{align*}
\|\Lambda^{r-\beta}(u\cdot\nabla\tau)\|_{L^2}
&\lesssim
\|u\|_{L^\infty}
\|\Lambda^{r-\beta+1}\tau\|_{L^2}
+\|\tau\|_{L^\infty}
\|\Lambda^{r-\beta+1}u\|_{L^2}\notag\\
&\lesssim\delta D_r^{\frac{1}{2}}.
\end{align*}
Consequently,
\[
|\mathcal{J}_2|
\le
\|\Lambda^{r-\beta}(u\cdot\nabla\tau)\|_{L^2}
\|\Lambda^{r+\beta}\tau\|_{L^2}
\lesssim\delta D_r.
\]
For $J_3$, using Lemma \ref{lemma6} and the Sobolev embeddings, we have
\begin{align*}
\|\Lambda^{r-\beta}Q(\nabla u,\tau)\|_{L^2}
&\lesssim
\|\tau\|_{L^\infty}
\|\Lambda^{r-\beta+1}u\|_{L^2}
+\|\nabla u\|_{L^{\frac3\beta}}
\|\Lambda^{r-\beta}\tau\|_{L^{\frac6{3-2\beta}}}
\notag\\
&\lesssim
\delta D_r^{\frac{1}{2}}
+\|\Lambda^{\frac52-\beta}u\|_{L^2}
 E_r^{\frac{1}{2}}.
\end{align*}
It follows from Young's inequality that
\begin{align*}
|\mathcal{J}_3|
&\le
\|\Lambda^{r-\beta}Q(\nabla u,\tau)\|_{L^2}
\|\Lambda^{r+\beta}\tau\|_{L^2} \le
\left(C\delta+\frac1{100}\right) D_r
+C\|\Lambda^{\frac52-\beta}u\|_{L^2}^2 E_r.
\end{align*}
Applying $\Lambda^{r-\beta}$ to the equations for $u$ and
$\operatorname{div}\tau$, and taking the corresponding
$L^2$ inner products, we obtain
\begin{align*}
&\frac{d}{dt}
k\langle\Lambda^{r-\beta}u,
        \Lambda^{r-\beta}\operatorname{div}\tau\rangle
+\frac{k}{2}\|\Lambda^{r-\beta+1}u\|_{L^2}^2\\
=&-k\langle\Lambda^{r-\beta}\mathbb P(u\cdot\nabla u),
              \Lambda^{r-\beta}\operatorname{div}\tau\rangle
-k\langle\Lambda^{r-\beta}\operatorname{div}
              (u\cdot\nabla\tau),\Lambda^{r-\beta}u\rangle\\
&-k\langle\Lambda^{r-\beta}\operatorname{div}
              Q(\nabla u,\tau),\Lambda^{r-\beta}u\rangle
-k\langle\Lambda^{r+\beta}\operatorname{div}\tau,
              \Lambda^{r-\beta}u\rangle
+k\|\Lambda^{r-\beta}\mathbb P\operatorname{div}\tau\|_{L^2}^2
\notag\\
\triangleq&\sum_{i=4}^{8}\mathcal{J}_i.
\end{align*}
Using $\operatorname{div}u=0$ and Lemma \ref{lemma6}, we find
\[
|\mathcal{J}_4|
\lesssim
\|u\|_{L^\infty}
\|\Lambda^{r-\beta+1}u\|_{L^2}\|\Lambda^{r-\beta+1}\tau\|_{L^2}
\lesssim\delta D_r.
\]
Similar to $|\mathcal{J}_2|$ and $|\mathcal{J}_3|$, we have $|\mathcal{J}_5|\lesssim\delta D_r$ and
\begin{align*}
|\mathcal{J}_6|
\le
\left(C\delta+\frac1{100}\right) D_r
+C \|\Lambda^{\frac52-\beta}u\|_{L^2}^2 E_r.
\end{align*}
Finally, using Young's inequality, we get
\begin{align*}
|\mathcal{J}_7|+|\mathcal{J}_8|
\le \frac{1}{100}\|\Lambda^{r+\beta}\tau\|_{L^2}^2+Ck^2\|\Lambda^{r-\beta+1} u\|_{L^2}^2+Ck\|\Lambda^{r-\beta+1}\tau\|_{L^2}^2.
\end{align*}
Combining the estimates above and choosing $k$ sufficiently small, we arrive at
\begin{align*}
    \frac{1}{2}\frac{d}{dt} E_r +D_r \le C\delta D_r+C\left(\delta^2+\|\Lambda^{\frac52-\beta}u\|_{L^2}^2\right) E_r.
\end{align*}
Noting that $\|\Lambda^{\frac52-\beta}u\|_{L^2}\lesssim\|\nabla u\|_{H^{s-\beta}}$, and choosing $\delta$ sufficiently small, we obtain from Gronwall's inequality that
\begin{align*}
    \sup_{t\in[0,T]} E_r(t)+\int_0^T D_r(t)dt\le C E_r(0).
\end{align*}
\end{proof}
Now we turn to the proof of continuous dependence.
Assume that $(u,\tau)$ and $(v,\sigma)$ are two solutions of (\ref{eq2}) with initial data $(u_0,\tau_0)\in H^s$ and $(v_0,\sigma_0)\in H^s$. Assume further that $(u_N,\tau_N)$ and $(v_N,\sigma_N)$ are two solutions of (\ref{eq2}) with initial data $(S_N u_0,S_N \tau_0)$ and $(S_N v_0,S_N \sigma_0)$. Then $\bar{u}\triangleq u-u_N,\ \bar{\tau}\triangleq \tau-\tau_N$ satisfy the following equations
\begin{align}
\left\{\begin{array}{l}
\partial_t\bar{u}+u\cdot\nabla \bar{u}+\bar{u} \cdot\nabla u_N + \nabla (P(u)-P(u_N)) =\rm{div}~\bar{\tau},~~~~{\rm div}~\bar{u}=0,\\[1ex]
\partial_t\bar{\tau}+u\cdot \nabla \bar{\tau}+\bar{u}\cdot \nabla\tau_N+Q(\nabla\bar{u},\tau)+Q(\nabla u_N,\bar{\tau})+\Lambda^{2\beta}\bar{\tau}=D(\bar{u}),\\[1ex]
\bar{u}|_{t=0}=u_0-S_N u_0,~~\bar{\tau}|_{t=0}=\tau_0-S_N \tau_0. \\[1ex]
\end{array}\right.
\end{align}
As before, we have 
\begin{align}
    \frac{1}{2}\frac{d}{dt}\left(\|\bar{u}\|_{H^s}^2+\|\bar{\tau}\|_{H^s}^2+2k\langle\bar{u},\rm{div} \bar{\tau}\rangle_{H^{s-\beta}}\right)+\frac{k}{2}\|\nabla\bar{u}\|_{H^{s-\beta}}^2+\|\Lambda^{\beta}\bar{\tau}\|_{H^{s}}^2
    \triangleq\sum_{i=1}^{14}L_i,
\end{align}
where 
\begin{align*}
    L_1&=-\langle\bar{u}\cdot \nabla u_N,\bar{u}\rangle-\langle\bar{u}\cdot \nabla \tau_N,\bar{\tau}\rangle
    -\langle Q(\nabla\bar{u},\tau),\bar{\tau}\rangle-\langle Q(\nabla u_N,\bar{\tau}),\bar{\tau}\rangle\\
    &~~~~~-k\langle \mathbb{P}(u \cdot \nabla\bar{u}),\mathrm{div}\bar{\tau}\rangle-k\langle \mathbb{P}(\bar{u}\cdot \nabla u_N),\mathrm{div}\bar{\tau}\rangle
    -k\langle \rm{div} (u \cdot \nabla \bar{\tau}),\bar{u}\rangle\\
    &~~~~~-k\langle \rm{div}(\bar{u}\cdot \nabla \tau_N),\bar{u}\rangle
    -k\langle\mathrm{div}Q(\tau,\nabla\bar{u}),\bar{u}\rangle-k\langle\mathrm{div}Q(\bar{\tau},\nabla u_N),\bar{u}\rangle,\\
    L_2&=-\langle\Lambda^s(u\cdot \nabla\bar{u}),\Lambda^s\bar{u}\rangle,~~~~~~~~~~~~~~~~~
    L_3=-\langle\Lambda^s(\bar{u}\cdot\nabla u_N),\Lambda^s \bar{u}\rangle,\\
    L_4&=-\langle\Lambda^s(u\cdot \nabla\bar{\tau}),\Lambda^s\bar{\tau}\rangle,
    ~~~~~~~~~~~~~~~~~L_{5}=-\langle\Lambda^s(\bar{u}\cdot\nabla \tau_N),\Lambda^s \bar{\tau}\rangle,\\
    L_{6}&=-\langle \Lambda^s Q(\nabla\bar{u},\tau),\Lambda^s\bar{\tau}\rangle,~~~~~~~~~~~~~~~L_{7}=-\langle \Lambda^s Q(\nabla u_N,\bar{\tau}),\Lambda^s\bar{\tau}\rangle,\\
    L_{8}&=-k\langle \Lambda^{s-\beta}\mathbb{P}(u\cdot\nabla\bar{u}),\Lambda^{s-\beta}\mathrm{div}\bar{\tau}\rangle,
    ~~~~~~L_{9}=-k\langle\Lambda^{s-\beta}\mathbb{P}(\bar{u}\cdot\nabla u_N),\Lambda^{s-\beta}\rm{div} \bar{\tau}\rangle,\\
    L_{10}&=-k\langle\Lambda^{s-\beta}\mathrm{div}(u\cdot\nabla \bar{\tau}),\Lambda^{s-\beta}\bar{u}\rangle,
    ~~~~~~~L_{11}=-k\langle\Lambda^{s-\beta}\rm{div}(\bar{u}\cdot \nabla \tau_N),\Lambda^{s-\beta}\bar{u}\rangle,\\
    L_{12}&=-k\langle\Lambda^{s-\beta}\mathrm{div}Q(\nabla \bar{u},\tau),\Lambda^{s-\beta}\bar{u}\rangle,
    ~~~~~L_{13}=-k\langle\Lambda^{s-\beta}\rm{div}Q(\nabla u_N,\bar{\tau}),\Lambda^{s-\beta}\bar{u}\rangle,\\
    L_{14}&=-k\langle \Lambda^{2\beta}\rm{div}\bar{\tau},\bar{u}\rangle-k\langle\Lambda^{s+\beta}\mathrm{div}\bar{\tau},\Lambda^{s-\beta}\bar{u}\rangle+k\langle \mathbb{P}\mathrm{div}\bar{\tau},\mathrm{div}\bar{\tau}\rangle+k\langle \Lambda^{s-\beta}\mathbb{P}\mathrm{div}\bar{\tau},\Lambda^{s-\beta}\mathrm{div}\bar{\tau}\rangle.
\end{align*}
It is easy to check that
\begin{align*}
    |L_1|\lesssim \left(\|\bar{u}\|_{H^s}^2+\|\bar{\tau}\|_{H^s}^2 \right)
    \left(\|u\|_{H^s}+\|\tau\|_{H^s}+\|u_N\|_{H^s}+\|\tau_N\|_{H^s} \right).
\end{align*}
We deduce from Lemma \ref{lemma4} that
\begin{align*}
    |L_2|&\lesssim \|\Lambda^s u\|_{L^{\frac{6}{2\beta+1}}}\|\nabla\bar{u}\|_{L^{\frac{3}{2-2\beta}}}\|\Lambda^s\bar{u}\|_{L^{\frac{6}{2\beta+1}}}+\|\nabla u\|_{L^{\frac{3}{2-2\beta}}}\|\Lambda^s\bar{u}\|_{L^{\frac{6}{2\beta+1}}}^2\\
    &\le C\|\Lambda^{s-\beta+1} u\|_{L^2}^2 \|\Lambda^{2\beta+\frac{1}{2}}\bar{u}\|_{L^2}^2
    +\frac{k}{100}\|\Lambda^{s-\beta+1}\bar{u}\|_{L^2}^2
    +C \epsilon^{\frac{1}{2}}\|\Lambda^{s-\beta+1} \bar{u}\|_{L^2}^2.
\end{align*}
In order to estimate $L_3$, we observe that, due to the Bony decomposition, Lemma \ref{lemma2} and Lemma \ref{lemma3},
\begin{align*}
    &\|\dot{T}_{\nabla u_N}\bar{u}\|_{\dot{B}^{s+\beta-1}_{2,2}}+\|\dot{R}\left( \nabla u_N, \bar{u}\right)\|_{\dot{B}^{s+\beta-1}_{2,2}}\lesssim \|\nabla u_N\|_{\dot{B}^{\beta-1}_{\infty,\infty}}\|\bar{u}\|_{\dot{B}^s_{2,2}}\lesssim\|u_N\|_{\dot{H}^{\beta+\frac{3}{2}}}\|\bar{u}\|_{\dot{H}^s}.
\end{align*}
Since $s\ge 2\beta+\frac12$ and $\frac12<\beta<1$, we have
\[
H^{s-\beta}
\hookrightarrow H^{\beta+\frac12}
\hookrightarrow L^{\frac{3}{1-\beta}}
\hookrightarrow \dot B_{\infty,\infty}^{\beta-1}.
\]
Applying Lemma \ref{lemma2} with $\beta-1<0$, we obtain that
\[
\begin{aligned}
\|\dot T_{\bar u}\nabla u_N\|_{\dot H^{s+\beta-1}}
&\lesssim
\|\bar u\|_{\dot B_{\infty,\infty}^{\beta-1}}
\|\nabla u_N\|_{\dot H^s}
\lesssim
\|\bar u\|_{H^{s-\beta}}
\|\Lambda^{s+1}u_N\|_{L^2}.
\end{aligned}
\]
Thus, it follows that
\begin{align*}
    |L_3|
    \le C\|\Lambda^{\beta+\frac{3}{2}}u_N\|_{L^2}^2 \|\bar{u}\|_{H^s}^2+\frac{k}{100}\|\Lambda^{s-\beta+1}\bar{u}\|_{L^2}^2+C\|\Lambda^{s+1}u_N\|_{L^2}^2\|\bar{u}\|_{H^{s-\beta}}^2.
\end{align*}
Similar to $L_2$ and $L_3$, we have 
\begin{align*}
    |L_4|\le C\|\Lambda^{s-\beta+1}u\|_{L^2}^2\|\Lambda^{2\beta+\frac{1}{2}}\bar{\tau}\|_{L^2}^2+\frac{1}{100}\|\Lambda^{s-\beta+1}\bar{\tau}\|_{L^2}^2+C \epsilon^{\frac{1}{2}}\|\Lambda^{s-\beta+1}\bar{\tau}\|_{L^2}^2,
\end{align*}
and
\begin{align*}
    |L_{5}|\le C\|\Lambda^{\beta+\frac{3}{2}}\tau_N\|_{L^2}^2\|\bar{u}\|_{H^s}^2+\frac{1}{100}\|\Lambda^{s-\beta+1}\bar{\tau}\|_{L^2}^2+C\|\Lambda^{s+1}\tau_N\|_{L^2}^2\|\bar{u}\|_{H^{s-\beta}}^2.
\end{align*}
By Hölder's and Young's inequalities, we find that
\begin{align*}
    |L_{6}|&\lesssim \|\Lambda^{s-\beta}\tau\|_{L^{6}}\|\nabla \bar{u}\|_{L^3}\|\Lambda^{s+\beta}\bar{\tau}\|_{L^2}+\|\tau\|_{L^{\infty}}\|\Lambda^{s-\beta+1}\bar{u}\|_{L^2}\|\Lambda^{s+\beta}\bar{\tau}\|_{L^2}\\
    &\le C\|\Lambda^{s-\beta+1}\tau\|_{L^2}^2\|\Lambda^{\frac{3}{2}}\bar{u}\|_{L^2}^2+\frac{1}{100}\|\Lambda^{s+\beta}\bar{\tau}\|_{L^2}^2+C\epsilon^{\frac{1}{2}}\|\Lambda^{s-\beta+1}\bar{u}\|_{L^2}^2,
    \end{align*}
    and
    \begin{align*}
    |L_{7}|&\lesssim \|\Lambda^{s-\beta}\bar{\tau}\|_{L^{\frac{6}{3-2\beta}}}\|\nabla u_N\|_{L^{\frac{3}{\beta}}}\|\Lambda^{s+\beta}\bar{\tau}\|_{L^2}+\|\bar{\tau}\|_{L^{\infty}}\|\Lambda^{s-\beta+1}u_N\|_{L^2}\|\Lambda^{s+\beta}\bar{\tau}\|_{L^2}\\
    &\le C\|\Lambda^{s}\bar{\tau}\|_{L^2}^2\|\Lambda^{\frac{5}{2}-\beta}u_N\|_{L^2}^2+\frac{1}{100}\|\Lambda^{s+\beta}\bar{\tau}\|_{L^2}^2+C\|\bar{\tau}\|_{H^s}^2\|\Lambda^{s-\beta+1}u_N\|_{L^2}^2.
\end{align*}
Replacing $\tau$ by $u$ in the first estimate above and $\bar{\tau}$ by $\bar{u}$ in the second, we get 
\begin{align*}
    |L_{8}|\le C\|\Lambda^{s-\beta+1}u\|_{L^2}^2\|\Lambda^{\frac{3}{2}}\bar{u}\|_{L^2}^2+\frac{1}{100}\|\Lambda^{s-\beta+1}\bar{\tau}\|_{L^2}^2+C\epsilon^{\frac{1}{2}}\|\Lambda^{s-\beta+1}\bar{u}\|_{L^2}^2,
\end{align*}
and
\begin{align*}
    |L_{9}| \le C\|\Lambda^{s}\bar{u}\|_{L^2}^2\|\Lambda^{\frac{5}{2}-\beta}u_N\|_{L^2}^2+\frac{1}{100}\|\Lambda^{s-\beta+1}\bar{\tau}\|_{L^2}^2+C\|\bar{u}\|_{H^s}^2\|\Lambda^{s-\beta+1}u_N\|_{L^2}^2.
\end{align*}
Estimating $L_{10}-L_{13}$ in the same manner, we have that
\begin{align*}
    \sum_{i=10}^{13}|L_i| \le &
    C \left( \|\bar{u}\|_{H^s}^2+\|\bar{\tau}\|_{H^s}^2\right)
    \left( \|\Lambda^{\frac{5}{2}-\beta}u\|_{L^2}^2+\|\Lambda^{\frac{5}{2}-\beta}\tau\|_{L^2}^2+\|\Lambda^{\frac{5}{2}-\beta}u_N\|_{L^2}^2+\|\Lambda^{s-\beta+1}u_N\|_{L^2}^2\right)\\
    &+\frac{k}{100}\|\Lambda^{s-\beta+1}\bar{u}\|_{L^2}^2
    +C \epsilon^{\frac{1}{2}}\|\Lambda^{s-\beta+1}\bar{\tau}\|_{L^2}^2.
\end{align*}
At last, we can derive directly that 
\begin{align*}
    |L_{14}| \le k^2C\|\nabla \bar{u}\|_{H^{s-\beta}}^2+
    \left(\frac{1}{100}+kC\right)\|\Lambda^{\beta}\bar{\tau}\|_{H^{s}}^2.
\end{align*}
Denote 
\begin{align*}
   & E(t)\triangleq \|(u(t),\tau(t))\|_{H^s}^2+2k\langle u,\rm{div}\tau\rangle_{H^{s-\beta}}~,~~~D(t)\triangleq\frac{k}{2}\|\nabla u(t)\|_{H^{s-\beta}}^2+\|\Lambda^{\beta}\tau(t)\|_{H^{s}}^2;~~~\\
    &\bar{E}(t)\triangleq\|(\bar{u}(t),\bar{\tau}(t))\|_{H^s}^2+2k\langle \bar{u},\rm{div}\bar{\tau}\rangle_{H^{s-\beta}}~,~~~\bar{D}(t)\triangleq\frac{k}{2}\|\nabla \bar{u}(t)\|_{H^{s-\beta}}^2+\|\Lambda^{\beta}\bar{\tau}(t)\|_{H^{s}}^2;\\
    &E_N(t) \triangleq\|(u_N(t),\tau_N(t))\|_{H^s}^2+2k\langle u_N,\rm{div}\tau_N\rangle_{H^{s-\beta}}~,~~D_N(t)\triangleq\frac{k}{2}\|\nabla u_N(t)\|_{H^{s-\beta}}^2+\|\Lambda^{\beta}\tau_N(t)\|_{H^{s}}^2.
\end{align*}
Adding up the estimates for $L_1-L_{14}$ and choosing k and $\epsilon$ sufficiently small,  we conclude that
\begin{align}{\label{4-26}}
    &\frac{1}{2}\frac{d}{dt} \bar{E}(t)+\bar{D}(t) \\ \nonumber
    \le&
    C\bar{E}(t)\left(E^{\frac{1}{2}}(t)+E_N^{\frac{1}{2}}(t)+D(t)+D_N(t)  \right)+C \|\Lambda^{s+1}(u_N,\tau_N)\|_{L^2}^2 \|\bar{u}\|_{H^{s-\beta}}^2.
\end{align}
Setting $r=s+\beta$ in Lemma \ref{highregularity}, we get that
\begin{align}{\label{4-27}}
    \int_0^T\|\Lambda^{s+1}(u_N,\tau_N)\|_{L^2}^2 dt\le C\|(S_N u_0,S_N\tau_0)\|_{H^{s+\beta}}^2\le C 2^{2\beta N}\|(u_0,\tau_0)\|_{H^s}^2\le C2^{2\beta N}.
\end{align}
Substituting $(\bar{u},\bar{\tau})$ for $(\delta u,\delta \tau)$ in (\ref{4-23}) yields that 
\begin{align}{\label{4-28}}
    \sup_{t\in [0,T]} \|\bar{u}\|_{H^{s-\beta}}^2 \le C\|(\bar{u}(0),\bar{\tau}(0))\|_{H^{s-\beta}}^2 \le C2^{-2\beta N}\|(\bar{u}(0),\bar{\tau}(0))\|_{H^s}^2.
\end{align}
Integrating (\ref{4-26}) and using (\ref{4-27}) and (\ref{4-28}), we get that 
\begin{align*}
    \bar{E}(t)+\int_0^t \bar{D}(t')dt'\le C\bar{E}(0)+ C\int_0^t \bar{E}(t')G(t')dt'.
\end{align*}
where $G(t)\triangleq E^{\frac{1}{2}}(t)+E_N^{\frac{1}{2}}(t)+D(t)+D_N(t)$ is integrable on [0,T].\\
An application of Gronwall's inequality yields
\begin{align}{\label{4-29}}
    \sup_{t\in [0,T]}\bar{E}(t)+\int_0^T \bar{D}(t)dt \le C \bar{E}(0),
\end{align}
and the similar estimate holds for $(v,\sigma)$.\\
We now turn to $u_N-v_N$ and $\tau_N-\sigma_N$. Substituting $u_N-v_N$ for $\delta u$ and $\tau_N-\sigma_N$ for $\delta\tau$ in (\ref{4-23}) and using Lemma \ref{highregularity}, we deduce that 
\begin{align}{\label{4-30}}
    &\|u_N-v_N\|_{L_T^{\infty}(H^s)\cap L_T^2(H^{s-\beta+1})}+\|\tau_N-\sigma_N\|_{L_T^{\infty}(H^s)\cap L_T^2(H^{s+\beta})}\\ \nonumber
    \le& C \|u_N-v_N\|_{L_T^{\infty}(H^{s-\beta})\cap L_T^2(H^{s-2\beta+1})}^{\frac{1}{2}}  \|u_N-v_N\|_{L_T^{\infty}(H^{s+\beta})\cap L_T^2(H^{s+1})}^{\frac{1}{2}}\\ \nonumber
    &+C \|\tau_N-\sigma_N\|_{L_T^{\infty}(H^{s-\beta})\cap L_T^2(H^{s})}^{\frac{1}{2}}  \|\tau_N-\sigma_N\|_{L_T^{\infty}(H^{s+\beta})\cap L_T^2(H^{s+2\beta})}^{\frac{1}{2}}  \\ \nonumber
   \le & C\left(\|S_N u_0-S_N v_0\|_{H^{s-\beta}}^{\frac{1}{2}}+\|S_N\tau_0-S_N\sigma_0\|_{H^{s-\beta}}^{\frac{1}{2}}\right)\\ \nonumber   &\cdot\left(\|u_N(0)\|_{H^{s+\beta}}^{\frac{1}{2}}+\|v_N(0)\|_{H^{s+\beta}}^{\frac{1}{2}}+\|\tau_N(0)\|_{H^{s+\beta}}^{\frac{1}{2}}+\|\sigma_N(0)\|_{H^{s+\beta}}^{\frac{1}{2}}\right)\\\nonumber
    \le& C2^{\frac{\beta N}{2}} \|(u_0-v_0,\tau_0-\sigma_0)\|_{H^s}^{\frac{1}{2}}.
\end{align}
Combining (\ref{4-29}) and (\ref{4-30}), we have
\begin{align*}
    &\|u-v\|_{L_T^{\infty}(H^s)\cap L_T^2(H^{s-\beta+1})}+\|\tau-\sigma\|_{L_T^{\infty}(H^s)\cap L_T^2(H^{s+\beta})}\\
    \le& C \|(Id-S_N) (u_0,v_0)\|_{H^s}+ C \|(Id-S_N) (\tau_0,\sigma_0)\|_{H^s}+C2^{\frac{\beta N}{2}}(\|u_0-v_0\|_{H^s}^{\frac{1}{2}}+\|\tau_0-\sigma_0\|_{H^s}^{\frac{1}{2}}).
\end{align*}
For any $\eta>0$, there exists an $N>0$ such that $C \|(Id-S_N) (u_0,\tau_0)\|_{H^s} \le \frac{\eta}{4}$. Then, there exists a $\gamma(\eta,N,u_0,\tau_0)>0$ such that for any $\|u_0-v_0\|_{H^s}+\|\tau_0-\sigma_0\|_{H^s}\le \gamma$, there holds
\begin{align*}
    C \|(Id-S_N) (v_0,\sigma_0)\|_{H^s} &\le C \|(Id-S_N) (u_0,\tau_0)\|_{H^s}+C\|(u_0-v_0,\tau_0-\sigma_0)\|_{H^s}\\
    &\le\frac{\eta}{4}+C\gamma \le \frac{\eta}{2},
\end{align*}
and
\begin{align*}
    C2^{\frac{\beta N}{2}}(\|u_0-v_0\|_{H^s}^{\frac{1}{2}}+\|\tau_0-\sigma_0\|_{H^s}^{\frac{1}{2}}) \le \frac{\eta}{4}.
\end{align*}
Therefore, one can infer that 
\begin{align*}
    \|u-v\|_{L_T^{\infty}(H^s)\cap L_T^2(H^{s-\beta+1})}+\|\tau-\sigma\|_{L_T^{\infty}(H^s)\cap L_T^2(H^{s+\beta})} \le \eta.
\end{align*}
According to Step 1 to Step 4, we have proved the local well-posedness of (\ref{eq2}) in the sense of Hadamard. 
\\

\section{Proof of Theorem 1.1 (ii)}
We now show that the solution constructed above is actually global under the condition $\beta\in(\frac{1}{2},\frac{5}{6}]$.\\
Assume that $(u,\tau)$ is a local solution of the system (\ref{eq2}) with the initial data satisfying (\ref{cond1}). Step 1 to Step 4 yield for any $t\in [0,T]$, $(u,\tau)$ satisfy 
\begin{align*}
    E_0(t) \le 2M,~~~~and~~~~E_1(t) \le 2\epsilon(M).
\end{align*}
We recall that
\[
E_0\triangleq \|(u,\tau)\|_{H^\beta}^2 ~~~~,
~~~
E_1\triangleq
\|\Lambda^\beta(u,\tau)\|_{H^{s-\beta}}^2.
\]
and further define
\[
D_0 \triangleq\frac{k}{2}\|\nabla u\|_{L^2}^2+\|\Lambda^{\beta}\tau\|_{H^{\beta}}^2~~~~\text{and}~~~~
D_1\triangleq\frac{k}{2}\|\Lambda^{2\beta}u\|_{H^{s-3\beta+1}}^2+\|\Lambda^{2\beta}\tau\|_{H^{s-\beta}}^2.
\]
Now, we shall prove the following two statements.
\begin{align}{\label{3-32}}
    \frac{d}{dt} E_1(t)+D_1(t)\le  C E_1(t)^{\frac{1}{2}} D_1(t),  
\end{align}
and 
\begin{align}{\label{3-33}}
    \frac{d}{dt}E_0(t)+D_0(t) \le C E_1(t)^{\frac{1}{2}}D_0(t).
\end{align}
\textbf{1.~The proof of (\ref{3-32})}
\\
The $\beta$-order energy can be written as 
\begin{align*}
&\quad\frac{1}{2}\frac{d}{dt}\|(\Lambda^{\beta} u,\Lambda^{\beta} \tau)\|_{L^{2}}^{2}+\|\Lambda^{2\beta}\tau\|_{L^{2}}^{2}\\\nonumber
&=-\langle \Lambda^{\beta} (u\cdot \nabla u),\Lambda^{\beta} u \rangle-\langle \Lambda^{\beta} (u\cdot \nabla \tau),\Lambda^{\beta} \tau\rangle
-\langle \Lambda^{\beta} Q(\nabla u,\tau), \Lambda^{\beta}\tau\rangle.
\end{align*}
Note that $\frac{5}{2}-2\beta \ge \beta$ since $\beta \le \frac{5}{6}$, from which it follows that
\begin{align*}
    |\langle \Lambda^{\beta} (u\cdot \nabla u),\Lambda^{\beta} u \rangle| &\lesssim
    \|u\|_{L^{\frac{3}{2\beta-1}}}\|\nabla u\|_{L^{\frac{6}{5-4\beta}}}\|\Lambda^{2\beta}u\|_{L^2}\\
    &\lesssim \|\Lambda^{\frac{5}{2}-2\beta}u\|_{L^2}\|\Lambda^{2\beta}u\|_{L^2}^2.
\end{align*}
Similarly, we have 
\[
|\langle \Lambda^{\beta} (u\cdot \nabla \tau),\Lambda^{\beta} \tau\rangle|
+|\langle \Lambda^{\beta} Q(\nabla u,\tau), \Lambda^{\beta}\tau\rangle|
\lesssim \|\Lambda^{\frac{5}{2}-2\beta}(u,\tau)\|_{L^2}\|\Lambda^{2\beta}(u,\tau)\|_{L^2}^2.
\] 
Hence, we can deduce that 
\begin{align}{\label{3-34}}
    \frac{1}{2}\frac{d}{dt}\|(\Lambda^{\beta} u,\Lambda^{\beta} \tau)\|_{L^{2}}^{2}+\|\Lambda^{2\beta}\tau\|_{L^{2}}^{2} \lesssim E_1^{\frac{1}{2}}D_1.
\end{align}
Just like $J'_1-J'_3$, the $s$-order energy can be handled as 
\begin{align}{\label{3-35}}
&\quad\frac{1}{2}\frac{d}{dt}\|(\Lambda^{s} u,\Lambda^{s} \tau)\|_{L^{2}}^{2}+\|\Lambda^{s+\beta}\tau\|_{L^{2}}^{2}\\\nonumber
&=-\langle \Lambda^{s} (u\cdot \nabla u),\Lambda^{s} u \rangle-\langle \Lambda^{s} (u\cdot \nabla \tau),\Lambda^{s} \tau\rangle
-\langle \Lambda^{s} Q(\nabla u,\tau), \Lambda^{s}\tau\rangle\\ \nonumber
&\lesssim E_1^{\frac{1}{2}}D_1.
\end{align}
Now we consider the inner product for $u$ and $\rm{div} \tau$. According to (\ref{eq2}) and (\ref{3-15}), we can deduce that
\begin{align}{\label{3-36}}
  &\frac{d}{dt}\langle \Lambda^{2\beta-1}u,k\Lambda^{2\beta-1}\mathrm{div}\tau\rangle_{L^{2}}+\frac{k}{2}\|\Lambda^{2\beta} u\|_{L^{2}}^{2} 
  \\ \nonumber
  =&-k\langle\Lambda^{2\beta-1}(u\cdot\nabla u),\Lambda^{2\beta-1}\mathbb{P}\mathrm{div}\tau\rangle
  -k\langle\Lambda^{2\beta-1}\mathrm{div}(u\cdot\nabla\tau),\Lambda^{2\beta-1}u\rangle\\\nonumber
  &-k\langle\Lambda^{2\beta-1}\mathrm{div} Q(\nabla u,\tau),\Lambda^{2\beta-1}u\rangle
  -k\langle\Lambda^{2\beta-1}(-\Delta)^\beta\mathrm{div}\tau,\Lambda^{2\beta-1}u\rangle\\\nonumber
  &+k\langle\mathbb{P}\Lambda^{2\beta-1}\mathrm{div}\tau, \Lambda^{2\beta-1}\mathrm{div}\tau\rangle\\\nonumber
  \le & C\|\Lambda^{2\beta-1}(u,\tau)\|_{L^6}\|\nabla (u,\tau)\|_{L^3}\|\Lambda^{2\beta} (u,\tau)\|_{L^2}+C\|(u,\tau)\|_{L^{\infty}}\|\Lambda^{2\beta}(u,\tau)\|_{L^2}^2\\\nonumber
  &+\frac{1}{100}\|\Lambda^{4\beta-1}\tau\|_{L^2}^2
+Ck^2\|\Lambda^{2\beta}u\|_{L^2}^2+Ck\|\Lambda^{2\beta}\tau\|_{L^2}^2\\\nonumber
\le& C E_1^{\frac{1}{2}}D_1+\frac{1}{100}\|\Lambda^{4\beta-1}\tau\|_{L^2}^2+Ck^2\|\Lambda^{2\beta}u\|_{L^2}^2+Ck\|\Lambda^{2\beta}\tau\|_{L^2}^2.
\end{align}
Along the same line, replacing $2\beta-1$ by $s-\beta$, we find
\begin{align}{\label{3-37}}
  &\frac{d}{dt}\langle \Lambda^{s-\beta}u,k\Lambda^{s-\beta}\mathrm{div}\tau\rangle_{L^{2}}+\frac{k}{2}\|\Lambda^{s-\beta}\nabla u\|_{L^{2}}^{2} 
  \\ \nonumber
  =&-k\langle\Lambda^{s-\beta}(u\cdot\nabla u),\Lambda^{s-\beta}\mathbb{P}\mathrm{div}\tau\rangle
  -k\langle\Lambda^{s-\beta}\mathrm{div}(u\cdot\nabla\tau),\Lambda^{s-\beta}u\rangle\\ \nonumber
  &-k\langle\Lambda^{s-\beta}\mathrm{div} Q(\nabla u,\tau),\Lambda^{s-\beta}u\rangle
  -k\langle\Lambda^{s-\beta}(-\Delta)^\beta\mathrm{div}\tau,\Lambda^{s-\beta}u\rangle\\ \nonumber
  &+k\langle\mathbb{P}\Lambda^{s-\beta}\mathrm{div}\tau, \Lambda^{s-\beta}\mathrm{div}\tau\rangle\\ \nonumber
  \le & CE_1^{\frac{1}{2}}D_1+\frac{1}{100}\|\Lambda^{s+\beta} \tau\|_{L^2}^2+Ck^2\|\Lambda^{s-\beta+1}u\|_{L^2}^2+Ck\|\Lambda^{s-\beta+1}\tau\|_{L^2}^2.
\end{align}
Combining (\ref{3-34})-(\ref{3-37}) and choosing $k$ small enough, then we deduce (\ref{3-32}).

\textbf{2.~The proof of (\ref{3-33})}
\\
 Since $\tau$ is symmetric, we have 
$(\tau\Omega-\Omega\tau):\tau=0$ and $D(u)\tau+\tau D(u):\tau=2D(u):\tau^2.$
Thus, together with the definition of $D(u)$, we can write
\[
\langle Q(\nabla u,\tau),\tau\rangle
=-2b\sum_{i,j}\int_{\mathbb R^3}
\partial_j u_i\,(\tau^2)_{ij}\,dx.
\]
Let $R_j$ denote the Riesz transform. Using $\partial_j=\Lambda^{1-s} R_j \Lambda^{s}$, we get that
\begin{align*}
\langle\partial_j u_i,(\tau^2)_{ij}\rangle
&=\langle\Lambda^{1-s}R_j\Lambda^s u_i,(\tau^2)_{ij}\rangle=\langle R_j\Lambda^s u_i,
                  \Lambda^{1-s}(\tau^2)_{ij}\rangle
\end{align*}
Hence, by the  $L^p$ boundedness of $R_j$, together with Lemma \ref{lemma6}, we deduce that
\begin{align*}
|\langle Q(\nabla u,\tau), \tau\rangle|
&\lesssim \|\Lambda^{\beta} u\|_{L^2}\|\Lambda^{1-\beta} \tau\|_{L^{\frac{12}{5-2\beta}}}\|\tau\|_{L^{\frac{12}{1+2\beta}}}\\
&\lesssim \|\Lambda^{\beta}u\|_{L^2}\|\Lambda^{\frac{5-2\beta}{4}}\tau\|_{L^2}^2.
\end{align*}
Since $\beta \le \frac{5}{6}$ implies $\frac{5-2\beta}{4} \ge \beta$, this yields
\begin{align}\label{3-38}
\quad\frac{1}{2}\frac{d}{dt}\|(u,\tau)\|_{L^{2}}^{2}+\|\Lambda^{\beta}\tau\|_{L^{2}}^{2}\lesssim E_1^{\frac{1}{2}} D_0.
\end{align}
It is not difficult to check that
\begin{align}\label{3-39}
    &\frac{1}{2}\frac{d}{dt}\|(\Lambda^{\beta} u,\Lambda^{\beta} \tau)\|_{L^{2}}^{2}+\|\Lambda^{2\beta}\tau\|_{L^{2}}^{2}
    \lesssim \|(u,\tau)\|_{L^6}\|\nabla (u,\tau)\|_{L^2}\|\Lambda^{2\beta}(u,\tau)\|_{L^3}
    \lesssim E_1^{\frac{1}{2}}D_0,
\end{align}
and
\begin{align}\label{3-40}
  &\frac{d}{dt}\langle u,k\mathrm{div}\tau\rangle+\frac{k}{2}\|\nabla u\|_{L^{2}}^{2}\\\nonumber
  =&-k\langle(u\cdot\nabla u),\mathbb{P}\mathrm{div}\tau\rangle  
  -k\langle\mathrm{div}(u\cdot\nabla\tau),u\rangle
  -k\langle\mathrm{div} Q(\nabla u,\tau),u\rangle
  -k\langle(-\Delta)^\beta\mathrm{div}\tau,u\rangle \\\nonumber
  &+k\langle\mathbb{P}\mathrm{div}\tau, \mathrm{div}\tau\rangle\\\nonumber
  \le&\ C\|u\|_{L^{\infty}} \|\nabla u\|_{L^2} \|\nabla\tau\|_{L^2}
    +C\|u\|_{L^{\infty}} \|\nabla \tau\|_{L^2} \|\nabla u\|_{L^2}
    +C\|\nabla u\|_{L^2} \|\tau\|_{L^{\infty}}\|\nabla u\|_{L^2} \\\nonumber
  &+k\|\Lambda^{2\beta}\tau\|_{L^2} \|\nabla u\|_{L^2}
    +Ck\|\nabla\tau\|^2_{L^2}
    \\\nonumber
  \le&\ C E_1^{\frac{1}{2}}D_0
  +Ck^2\|\nabla u\|_{L^2}^2
  +\frac{1}{100}\|\Lambda^{2\beta}\tau\|_{L^2}^2
  +Ck\|\nabla\tau\|^2_{L^2}.
\end{align}
Combining (\ref{3-38}) to (\ref{3-40}) and choosing $k$ sufficiently small, we then get (\ref{3-33}).
\\
A bootstrap argument applied to (\ref{3-32}) yields that $\sup_{t\ge0} E_1(t) \le C\epsilon$ for some suitable $\epsilon>0$. Together with (\ref{3-33}), this gives us that if $\epsilon$ is sufficiently small, then for all $t>0$ we have
\begin{align*}
    E_0(t)+\int_0^t D_0(t')dt'\le C E_0(0)
\end{align*}
and the result follows.
\end{proof}

\section{Proof of Corollary \ref{coro1.2}}
	\par
In this section, we consider the critical case $s=2\beta+\frac{1}{2}$ with $\beta\in(\frac{1}{2},1)$ and small initial data.
\begin{proof}

Denote
\begin{align*}
    \mathcal{E}&\triangleq \|(u,\tau)\|_{H^{2\beta+\frac{1}{2}}}^2+2k\langle u, \rm{div} \tau\rangle_{H^{\beta+\frac{1}{2}}},\\
    \mathcal{D}&\triangleq \frac{k}{2}\|\nabla u\|_{H^{\beta+\frac{1}{2}}}^2+\|\Lambda^{\beta}\tau\|_{H^{2\beta+\frac{1}{2}}}^2.
\end{align*}
Note that 
\begin{align}{\label{5-1}}
\frac{1}{2}\frac{d}{dt}\|(u,\tau)\|_{L^{2}}^{2}+\|\Lambda^{\beta}\tau\|_{L^{2}}^{2}=-\langle Q(\nabla u,\tau), \tau\rangle\lesssim\|\nabla u\|_{L^2}\|\tau\|_{L^3}\|\tau\|_{L^6}\lesssim\mathcal{E}^{\frac{1}{2}} \mathcal{D}.
\end{align}
Substituting $s=2\beta+\frac{1}{2}$ into (\ref{3-35}) and (\ref{3-37}) and using the small initial data condition, together with (\ref{5-1}) and (\ref{3-40}), then we deduce that
\begin{align*}
    \frac{d}{dt}\mathcal{E}(t)+\mathcal{D}(t) \le C \mathcal{E}(t)^{\frac{1}{2}} \mathcal{D}(t).
\end{align*}
Since the initial data is small enough, the result follows from a standard bootstrap argument.
\end{proof}

\smallskip
\noindent\textbf{Acknowledgments}
This work was partially supported by the National Natural Science Foundation of
  China (No.12571261).

\vspace{2mm}

\noindent\textbf{Conflict of interest.} The authors do not have any possible conflicts of interest.

\vspace{2mm}

\noindent\textbf{Data availability statement.}
 Data sharing is not applicable to this article, as no data sets were generated or analyzed during the current study.

\phantomsection
\addcontentsline{toc}{section}{\refname}
\bibliographystyle{abbrv} 
\bibliography{Oldroyd-Bref}

\end{document}